%% file: article3.tex
\documentclass[11pt, reqno]{amsart}

\usepackage[T1]{fontenc}

\usepackage{natbib}

\usepackage{amsmath, amsfonts, amssymb, amsthm}

\usepackage{graphicx}

\usepackage{ae}
\usepackage[cm]{aeguill}

\input{notation.tex}

\usepackage[english]{babel}

\numberwithin{equation}{section}

\begin{document}

\title[Sharp estimation with random design]{Sharp estimation in sup
  norm with random design}

\author[S. Ga{\"\i}ffas]{ St{\'e}phane Ga{\"\i}ffas }

\address{Laboratoire de Probabilit{\'e}s et Mod{\`e}les
  Al{\'e}atoires \\
  Universit{\'e} Paris 7, 175 rue du Chevaleret, 75013 Paris \\
}

\email{gaiffas@math.jussieu.fr}



\date{\today}

\dedicatory{ Laboratoire de Probabilit{\'e}s et Mod{\`e}les
  Al{\'e}atoires \\
  Universit{\'e} Paris 7, 175 rue du Chevaleret, 75013 Paris \\
  email: \rm{ \texttt{gaiffas@math.jussieu.fr} } }

\setlength{\leftmargini}{1.3em}

\begin{abstract}
  The aim of this paper is to recover the regression function with sup
  norm loss. We construct an asymptotically sharp estimator which
  converges with the spatially dependent rate
  \begin{equation*}
    r_{n, \mu}(x) = P \big(\log n / (n \mu(x)) \big)^{s / (2s + 1)},
  \end{equation*}
  where $\mu$ is the design density, $s$ the regression smoothness,
  $n$ the sample size and $P$ is a constant expressed in terms of a
  solution to a problem of optimal recovery as in~\cite{donoho94}.  We
  prove this result under the assumption that $\mu$ is positive and
  continuous.  This estimator combines kernel and local polynomial
  methods, where the kernel is given by optimal recovery, which allows
  to prove the result up to the constants for any $s > 0$.  Moreover,
  the estimator does not depend on $\mu$.  We prove that $r_{n,
    \mu}(x)$ is optimal in a sense which is stronger than the
  classical minimax lower bound.  Then, an inhomogeneous confidence
  band is proposed.  This band has a non constant length which depends
  on the local amount of data.
\end{abstract}

\maketitle

\section{Introduction \& main results}
\label{sec:introduction}

\subsection{The model}

Suppose we observe $(X_i, Y_i), 1 \leq i \leq n$, from
\begin{equation}
  \label{eq:model_regression}
  Y_i = f(X_i) + \xi_i,
\end{equation}
where $\xi_i$ are i.i.d. centered Gaussian with variance $\sigma^2$
and independent of $X_i$, with $X_i$ i.i.d. with density $\mu$ on $[0,
1]$, which is bounded away from $0$. We want to recover~$f$. In this
model, when $\mu$ is not the uniform law, we say that the information
is spatially inhomogeneous.

\subsection{Methodology}

There are several ways to assess the quality of an estimation
procedure. A first approach is local: we focus on recovering $f$ at a
fixed point $x_0 \in [0, 1]$. Over a function class $\Sigma$, the
minimax risk is given by
\begin{equation*}
  \mc R_n(\Sigma, x_0) = \inf_{\wh f_n} \sup_{f \in \Sigma} \mb E_f^n
  \big\{ |\wh f_n(x_0) - f(x_0) | \big\},
\end{equation*}
where the infimum is taken among all estimators. We say that
$\rho_n(x_0) > 0$ is the minimax convergence rate at $x_0$ if
\begin{equation*}
  0 < \liminf_n \frac{\mc R_n(\Sigma, x_0)}{\rho_n(x_0)} \leq
  \limsup_n \frac{\mc R_n(\Sigma, x_0)}{\rho_n(x_0)} < +\infty.
\end{equation*}
In this paper, we are interested in recovering $f$ globally. We
consider the loss with sup norm defined by $\norminfty{g} = \sup_{x
  \in [0, 1]} |g(x)|$. In this case, the minimax risk is
\begin{equation}
  \label{eq:R_n_sup_chap3}
  \mc R_n(\Sigma) = \inf_{\wh f_n} \sup_{f \in \Sigma} \mb E_f^n
  \big\{\norminfty{\wh f_n - f} \big\},
\end{equation}
and we say that $\psi_n$ is the minimax convergence rate if
\begin{equation*}
  0 < \liminf_n \frac{\mc R_n(\Sigma)}{\psi_n} \leq \limsup_n
  \frac{\mc R_n(\Sigma)}{\psi_n} < +\infty.
\end{equation*}
An advantage of this norm is that it is exacting: it forces an
estimator to behave well at every point simultaneously. In the
regression model~\eqref{eq:model_regression} with $\Sigma$ a H\"older
ball with smoothness $s > 0$, we have when $\mu$ is positive and
bounded that $\psi_n \asymp (\log n / n)^{s / (2s + 1)}$ (see
\cite{stone82}), where $a_n \asymp b_n$ means $0 < \liminf_n a_n /
b_n\leq \limsup_n a_n / b_n < +\infty$.

However, when $\mu$ is positive and bounded, $\psi_n$ is not sensitive
to the variations in the amount of data. An improvement is to consider
instead of~\eqref{eq:R_n_sup_chap3} the spatially dependent risk
\begin{equation*}
  \sup_{f \in \Sigma} \mb E_f^n \big\{ \sup_{x \in [0, 1]} r_n(x)^{-1}
  |\wh f_n(x) - f(x) | \big\},
\end{equation*}
where $\wh f_n$ is some estimator and $r_n(\cdot) > 0$ a family of
spatially dependent normalisation factors. If this quantity is bounded
as $n$ goes to infinity, we say that $r_n(\cdot)$ is an upper bound
over $\Sigma$. If we look for such upper bounds, we clearly find that
$r_n(x) \asymp \psi_n$ for any $x$, thus we must sharp this upper
bound up to constants. Here, we consider indeed the latter approach in
the asymptotic minimax context.  In this paper, we develop the
consequences of inhomogeneous data within this framework.

\subsection{Upper and lower bounds}

If $s, L > 0$, we define the H\"older ball $\Sigma(s, L)$, which is
the set of all the functions $f : [0, 1] \raro \setR$ such that for
any $x, y \in [0, 1]$,
\begin{equation*}
  |f^{(k)}(x) - f^{(k)}(y) | \leq L |x - y|^{s - k},  
\end{equation*}
where $k = \ppint{s}$ is the largest integer $k < s$. If $Q > 0$, we
denote by $\Sigma^Q(s, L)$ the set of functions $f \in \Sigma(s, L)$
such that $\norminfty{f} \leq Q$, and we denote simply $\Sigma =
\Sigma^Q(s, L)$. All along this study, we suppose:%
\setcounter{assumption}{3}%
\begin{assumption}
  \label{ass:design}
  For some $0 < \nu \leq 1$ and $\varrho, q > 0$, we have
  \begin{equation*}
    \mu \in \Sigma(\nu, \varrho) \, \text{ and } \, \mu(x) \geq
    q, \text{ for all } x \in [0, 1].
  \end{equation*}
\end{assumption}

In the following, a loss function $w(\cdot)$ is any non negative and
nondecreasing function such that $w(x) \leq A( 1 + |x|^b )$ for some
$A, b > 0$ (an example is $w(\cdot) = |\cdot|^p$ for $p > 0$).  Let us
consider
\begin{equation}
  \label{eq:r_n_def}
  r_{n, \mu}(x) = \Big( \frac{\log n}{n \mu(x)} \Big)^{s / (2s + 1)}.
\end{equation}
We denote by $\Efm$ the integration with respect to the joint law
$\Pfm$ of the observations $(X_i, Y_i)$, $1 \leq i \leq n$. Our first
result shows that $r_{n, \mu}(\cdot)$ is, up to the constants, an
upper bound over $\Sigma$.

\begin{theorem}[Upper bound]
  \label{thm:upper_bound}
  Under assumption \ref{ass:design}, if $\wh f_n$ is the estimator
  defined in section \ref{sec:estimator}, we have for any $s, L > 0$,
  \begin{equation}
    \label{eq:main_upper_bound}
    \limsup_n \sup_{f \in \Sigma} \Efm \big\{ w \big( \sup_{x \in [0,
      1] } r_{n, \mu}(x)^{-1}| \wh f_n(x) - f(x) | \, \big) \big\}
    \leq w(P),
  \end{equation}
  where
  \begin{equation}
    \label{eq:P_s_definition}
    P = \sigma^{2s / (2s + 1)} L^{1 / (2s + 1)} \, \varphi_s(0) \Big(
    \frac{2}{2s+1} \Big)^{s/(2s+1)}
  \end{equation}
  and $\varphi_s$ is defined as the solution of the optimisation
  problem
  \begin{equation}
    \label{eq:OR_optimal_function}
    \varphi_s \eqdef \argmax_{ \substack{ \varphi \in \Sigma(s, 1;
        \setR), \\ \norm{\varphi}_2 \leq 1 } } \varphi(0),
  \end{equation}
  where $\Sigma(s, L; \setR)$ is the extension of $\Sigma(s, L)$ to
  the whole real line.
\end{theorem}

In the same fashion as in~\cite{donoho94}, the constant $P$ is defined
via the solution of an optimisation problem which is connected to
optimal recovery. For further details, see in
sections~\ref{sec:discussion} and~\ref{sec:optimal_recovery}. The next
theorem shows that $r_{n, \mu}(\cdot)$ is indeed optimal in an
appropriate sense. In what follows, the notation $|I|$ stands for the
length of an interval $I$.

\begin{theorem}[Lower bound]
  \label{thm:lower_bound}
  Under assumption \ref{ass:design}, if $I_n \subset [0, 1]$ is any
  interval such that for some $\von \in (0, 1)$,
  \begin{equation}
    \label{eq:I_n_hyp}
    |I_n| \, n^{\von / (2s + 1)} \raro +\infty \quad \text{ as }\, n \raro
    +\infty,
  \end{equation}
  we have
  \begin{equation*}
    \liminf_n\, \inf_{\wh f_n} \sup_{f \in \Sigma} \Efm \big\{ w \big(
    \sup_{x \in I_n} r_{n, \mu}(x)^{-1} |\wh f_n(x) - f(x)| \big)
    \big\} \geq w\big((1 - \von) P\big),
  \end{equation*}
  where $P$ is given by \eqref{eq:P_s_definition} and the infimum is
  taken among all estimators. A consequence is that if $I_n$ is such
  that \eqref{eq:I_n_hyp} holds for any $\von \in (0, 1)$, we have
  \begin{equation}
    \label{eq:stronger_lower_bound}
    \liminf_n\, \inf_{\wh f_n} \sup_{f \in \Sigma} \Efm \big\{
    w \big( \sup_{x \in I_n} r_{n, \mu}(x)^{-1} |\wh f_n(x) - f(x)|
    \big) \big\} \geq w(P).
  \end{equation}
\end{theorem}

This result is discussed in details in
section~\ref{sec:about_the_lower_bound}. Now, we construct a
confidence band which is adapted to inhomogeneous data. Indeed, its
length varies depending on the local amount of data.

\subsection{An inhomogeneous confidence band}
\label{sec:confidence_bands}

We define the empirical design sample distribution
\begin{equation*}
  \bar \mu_n = \frac{1}{n} \sumin \delta_{X_i},
\end{equation*}
where $\delta$ is the Dirac mass, and for $h > 0$, $x \in [0,1]$, we
consider the intervals
\begin{equation}
  \label{eq:I_x_h_def}
  I(x, h) =
  \begin{cases}
    \,\, [x, x + h] &\text{ when } 0 \leq x \leq 1/2, \\
    \,\, [x - h, x] &\text{ when } 1/2 < x \leq 1.
  \end{cases}
\end{equation}
The choice of non symmetrical intervals allows to skip boundaries
effects. Then, we define the "bandwidth" at $x$ by
\begin{equation}
  \label{eq:data_driven_bandwidth_def}
  H_{n}(x) \eqdef \argmin_{h \in [0, 1]} \bigg \{ h^s \geq \Big(
  \frac{\log n}{n \bar \mu_n\big( I(x, h) \big)} \Big)^{1/2} \bigg \},
\end{equation}
which makes the balance between the bias and the variance of a certain
kernel estimator (more in section~\ref{sec:estimator} below). We
consider the sequence of points
\begin{equation}
  \label{eq:discretisation_points}
  x_j = j \Delta_n, \quad \Delta_n = (\log n)^{-2s / (2s + 1)} n^{-1 /
    (2s + 1)}, 
\end{equation}
for $j \in \mc J_n \eqdef \{ 0, \ldots, [\Delta_n^{-1}] \}$ where
$[a]$ is the integer part of $a$ with $x_{M_n} = 1$, $M_n = |\mc J_n|$
(the notation $|A|$ stands also for the size of a finite set $A$). If
$x \in [x_j, x_{j+1})$, we define
\begin{equation*}
  R_n(x) = H_n(x_j)^s,
\end{equation*}
and for any $x \in [0, 1]$, $\beta > 0$, we consider the band

\begin{equation}
  \label{eq:confidence_band_def}
  C_{n, \beta}(x) = \big[ \wh f_n(x) - (1 + \beta) P \, R_n(x),\,\,
  \wh f_n(x) + (1 + \beta) P \, R_n(x) \big],
\end{equation}
where $P$ is defined by~\eqref{eq:P_s_definition}. The next
proposition provides a control over the coverage probability of this
band, uniformly over $[0, 1]$.

\begin{proposition}
  \label{prop:CI}
  Given a confidence level $\alpha \in (0, 1)$, $C_{n, \beta}$ with
  \begin{equation*}
    \beta = \beta(n, \alpha) = \Big( \frac{\log(1 / \alpha)}{D_c (\log
      n)^{2s / (2s + 1)}} \Big)^{1/2}
  \end{equation*}
  \tup(where $D_c$ is some positive constant\tup), is under assumption
  \ref{ass:design}, a confidence band of asympotic level $1 - \alpha$,
  namely\tup:
  \begin{equation}
    \label{eq:CI_coverage_prob}
    \inf_{f \in \Sigma} \Pfm\big\{ \, f(x) \in C_{n, \beta}(x), \,
    \text{ for all } \, x \in [0, 1] \, \big\} \geq 1 - \alpha,
  \end{equation}
  for $n$ large enough. Moreover, we have for any $x \in [0, 1]$,
  \begin{equation}
    \label{eq:CI_length_equivalent}
    \sup_{f \in \Sigma} \Efm\{ | C_{n, \beta}(x) | \} / r_{n, \mu}(x)
    \raro 2 P\, \text{ as } n \raro +\infty.
  \end{equation}
\end{proposition}

In figures~\ref{fig_CI_homo} and~\ref{fig_CI_inhomo}, we give a
numerical illustration of this confidence band. We consider the
function $f(x) = 0.3 (1 - |x - 0.5| / 0.3)_+$, where $a_+ = \max(a,
0)$. The first dataset is simulated with an uniform design and the
second dataset with design density $\mu(x) = 0.05 + 11.4 |x -
0.5|^{2}$. In this example $s = L = 1$, the sample size is $n = 500$
and the root-signal-to-noise ratio is $7$.

\newlength{\figurelength}
\setlength{\figurelength}{6.5cm}

\begin{figure}[htbp]
  \begin{center}
    \includegraphics[width =
    \figurelength]{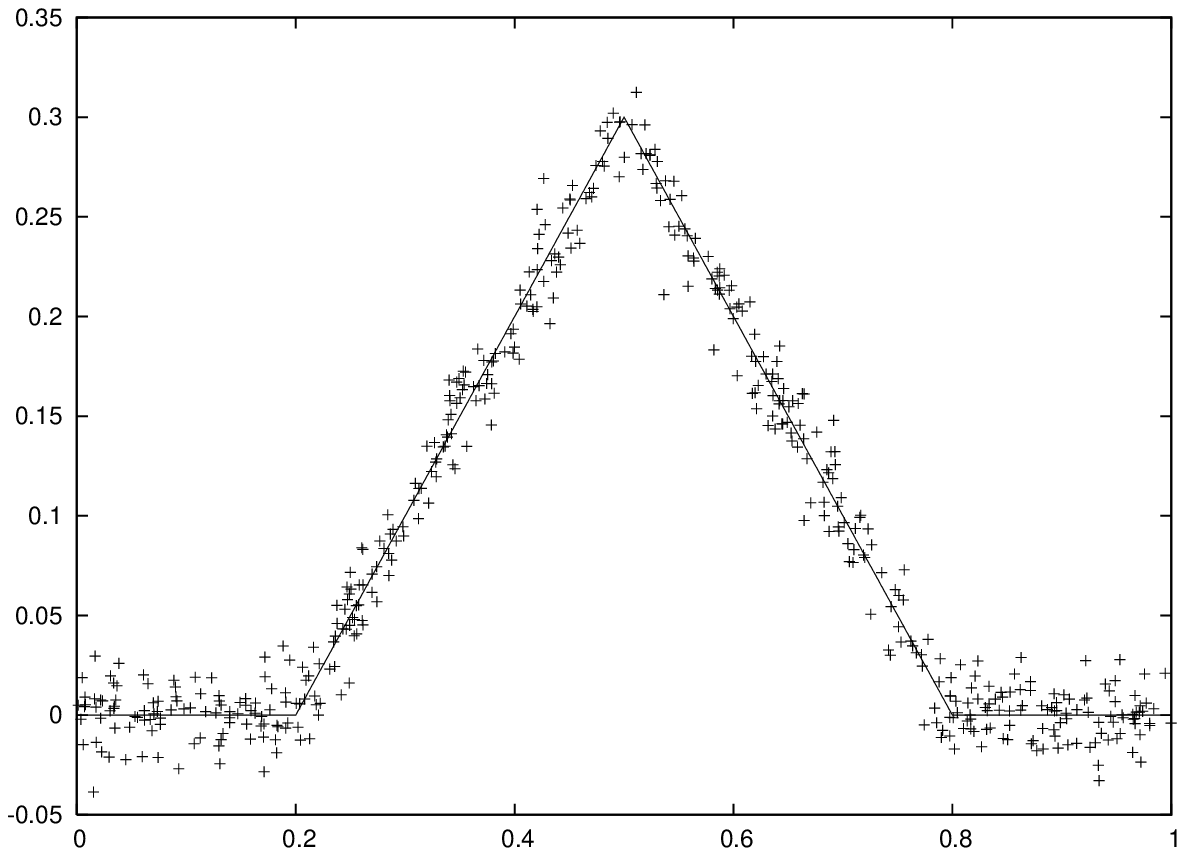}%
    \includegraphics[width =
    \figurelength]{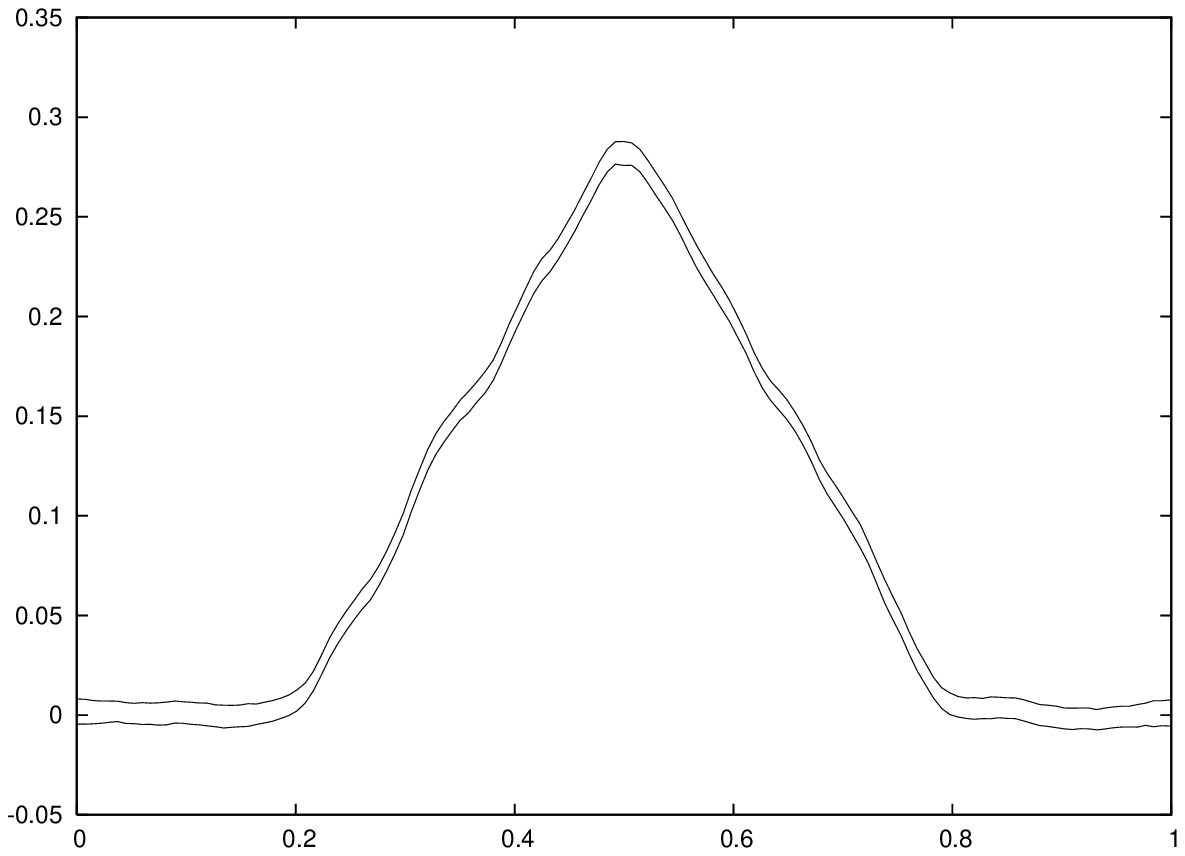}%
  \end{center}
  \caption{Confidence band with homogeneous data.}
  \label{fig_CI_homo}
\end{figure}

\begin{figure}[htbp]
  \begin{center}
    \includegraphics[width =
    \figurelength]{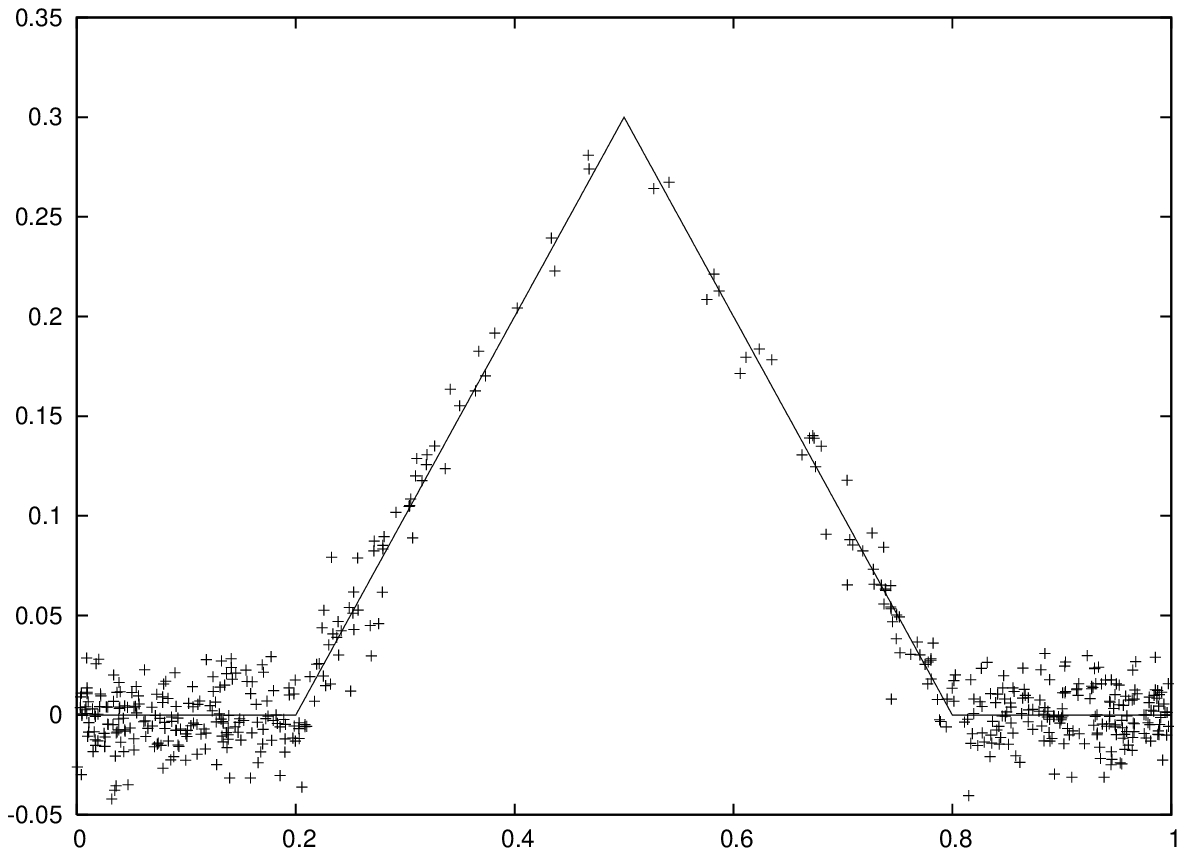}%
    \includegraphics[width =
    \figurelength]{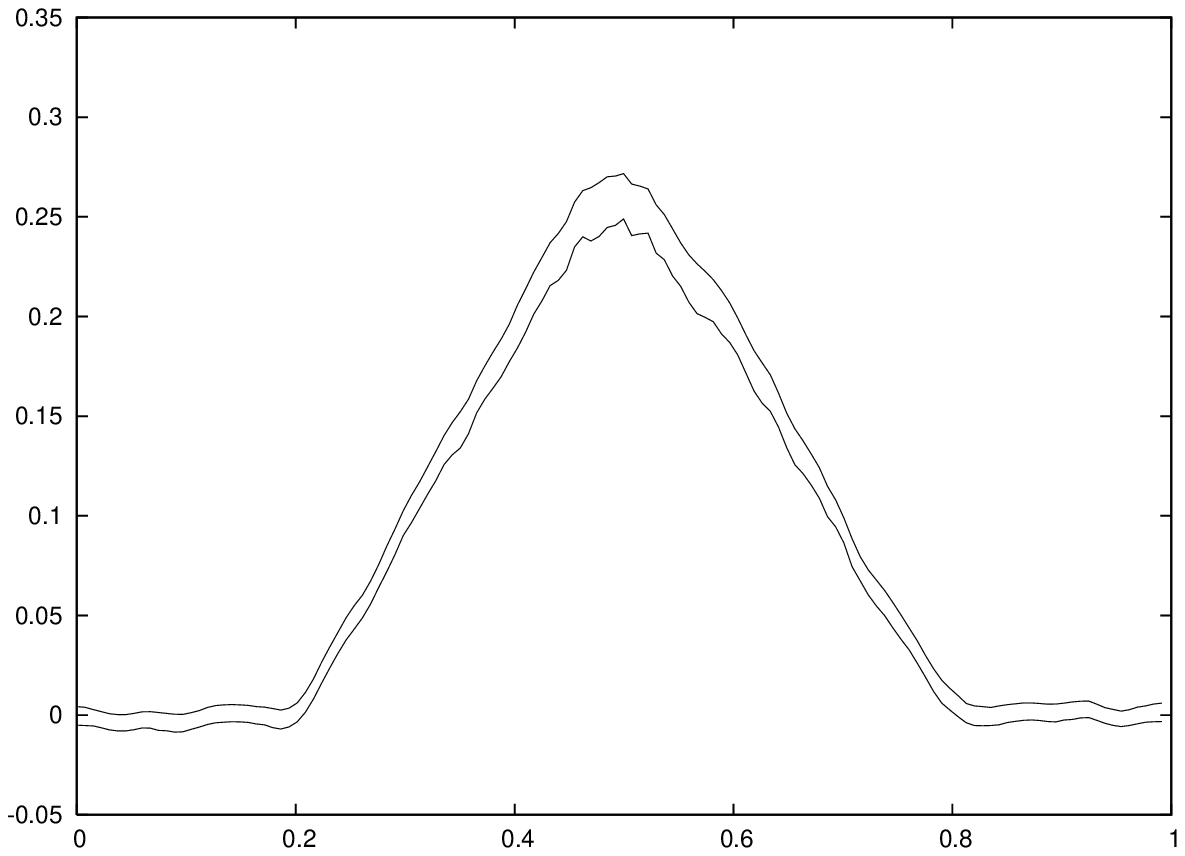}
  \end{center}
  \caption{Confidence band with inhomogeneous data.}
  \label{fig_CI_inhomo}
\end{figure}

When the data is homogeneous (uniform design), the length of the
confidence band is almost constant, see figure~\ref{fig_CI_homo}. In
the non-uniform case, the band is confined at the boundaries of $[0,
1]$ and more spaced at the middle, see figure~\ref{fig_CI_inhomo}.



\subsection{Outline}

The remainder of the paper is organised as follows. In
section~\ref{sec:discussion} we discuss our results in details and
compare them with former results. In section~\ref{sec:estimator}, we
construct the estimator used in theorem~\ref{thm:upper_bound}. The
proofs are delayed until sections~\ref{sec:proof_upper_bound}
and~\ref{sec:proof_of_the_lower_bound}. In
section~\ref{sec:optimal_recovery}, we recall some well known facts on
optimal recovery, which are useful for the construction of our
estimator and for the proofs.

\section{Discussion}
\label{sec:discussion}

\subsection{Motivation}

In most cases, the models considered in curve estimation do not allow
situations where the data is inhomogeneous, in so far as the amount of
data is implicitly assumed constant over space (or time).  However, an
increasing literature works in models where the data can be
inhomogeneously distributed. Recent results deal with the estimation
of the regression function when the observation points are not
equispaced or random, see for instance~\cite{antoniadis_et_al97},
\cite{cai_brown98}, \cite{wong_zheng02}, \cite{voichitaphd}, among
others. The estimators proposed in these papers present good minimax
properties, but the results are always stated in a way in which the
following basic principle does not appear: \emph{an estimator shall
  behave better at a point where there is much data than where there
  is little data}. For instance, upper bounds are usually stated with
the minimax rate, which is not sensitive to the variations in the
local amount of data nor to the information distribution in the
considered model.

At this stage, it is also natural to look for confidence bands when
the data is inhomogeneous, and especially distributed with an unknown
density. Following the above principle, a striking question is that of
the construction of a confidence band with a length which depends on
the local amount of data: such a band should be more confined where
there is much data than where there is little data. The aim of this
paper is to develop this new approach.


\subsection{Literature}

When the design is equidistant, that is $X_i = i/n$, we know
from~\cite{korostelev93} the exact asymptotic value of the minimax
risk for sup norm error loss. If
\begin{equation*}
  \psi_n = \Big( \frac{\log n}{n} \Big)^{s / (2s + 1)},
\end{equation*}
we have for any $0 < s \leq 1$ and $\Sigma = \Sigma(s, L)$,
\begin{equation*}
  \lim_{n \raro +\infty} \inf_{\wh f_n} \sup_{f \in \Sigma} \mbb E_f
  \big\{ w( \psi_{n}^{-1} \norminfty{\wh f_n - f} ) \big\} = w(C),
\end{equation*}
where
\begin{equation}
  \label{eq:C_1_def}
  C = \sigma^{2s / (2s + 1)} L^{1 / (2s + 1)} \Big( \frac{s + 1}{2
    s^2} \Big)^{s / (2s +  1)}.
\end{equation}
This result was the first of its kind for sup norm error loss. The
exact asymptotic value of the minimax risk was only known for square
integrated norm error loss, see~\cite{pinsker80}.

In the white noise model
\begin{equation}
  \label{eq:white_noise_model}
  d Y_t^{n} = f(t) dt + n^{-1/2} d W_t, \quad t \in [0,1],
\end{equation}
where $W$ is a standard Brownian motion,~\cite{donoho94} extends the
result by~\cite{korostelev93} to any $s > 1$. In this paper, the
author makes a link between statistical sup norm estimation and the
theory of optimal recovery (see section~\ref{sec:optimal_recovery}).
It is shown for any $s > 0$ and $\Sigma = \Sigma(s, L)$ that the
minimax risk satisfies
\begin{equation}
  \label{eq:donoho_result}
  \lim_{n \raro +\infty} \inf_{\wh f_n} \sup_{f \in \Sigma} \mb E_f
  \big\{ \psi_n^{-1} \norminfty{\wh f_n - f} \big\} = w(P_1),
\end{equation}
where $P_1$ is given by~\eqref{eq:P_s_definition} with $\sigma = 1$.
When $s \in (0, 1]$, we have $P = C$, see for instance
in~\cite{leonov97}.

Since the results by Korostelev and Donoho, many other authors worked
on the problem of sharp estimation (or testing) in sup norm. On
testing, see~\cite{lepski_tsybakov00},
see~\cite{korostelev_nussbaum_99} for density estimation
and~\cite{bertin04} for white noise in an anisotropic setting.

While most papers on sharp estimation work in models with homogeneous
information, the paper by~\cite{bertin02} works in the model of
regression with random design~\eqref{eq:model_regression}. When $\mu$
satisfies assumption~\ref{ass:design} and $\Sigma = \Sigma^Q(s, L)$
for $0 < s \leq 1$, it is shown that
\begin{equation}
  \label{eq:KB_minimax}
  \lim_{n \raro +\infty} \inf_{\wh f_n} \sup_{f \in \Sigma} \Efm
  \big\{ w( v_{n, \mu}^{-1} \norminfty{\wh f_n - f}) \big\} = w(C),
\end{equation}
where $C$ is given by~\eqref{eq:C_1_def} and
\begin{equation}
  \label{eq:KB_rate_def}
  v_{n, \mu} = \Big( \frac{\log n}{n \inf_x \mu(x)} \Big)^{s / (2s +
    1)}.
\end{equation}
Note that the rate $v_{n, \mu}$ differs from (and is larger than)
$\psi_n$ when $\mu$ is not uniform. A disappointing fact is that
$v_{n, \mu}$ depends on $\mu$ only via its infimum, which corresponds
to the point in $[0, 1]$ where we have the least information. This
rate does not take into account the regions with more data. It seems
natural to wonder if we can improve this result, namely: \emph{can we
  replace $\inf \mu$ by $\mu(x)$} ? Note that in
section~\ref{sec:introduction}, we have answered positively to this
question.

In this paper, we extend the result by~\cite{donoho94} to the model of
regression with random design and we improve the result
by~\cite{bertin02} in several ways: our result holds for any $s > 0$,
we construct an estimator which does not depend on $\mu$, and when the
design is not uniform, our convergence rate $r_{n, \mu}(\cdot)$ is
better (smaller) than $v_{n, \mu}$ at the order of constants.  More
importantly, this rate is adapted to the local amount of information
of the model.

\subsection{About theorem \ref{thm:upper_bound}}

We can understand the result of theorem~\ref{thm:upper_bound}
heuristically. Following \cite{brown_low96} and
\cite{brown_cai_low_zhang02} we can find an "idealised" statistical
experiment which is equivalent (in the sense that the LeCam deficiency
goes to $0$) to the model~\eqref{eq:model_regression}. The
model~\eqref{eq:model_regression} is clearly equivalent to
\begin{equation*}
  Y_i = f(G_{\mu}^{-1}(U_i)) + \xi_i, \quad 1 \leq i \leq n,
\end{equation*}
with independent and uniform $U_i$ where $G_{\mu}(x) = \int_0^x \mu(t)
dt$. Under appropriate conditions on $f$ and $\mu$, we know
from~\cite{brown_cai_low_zhang02} that this model is equivalent to
\begin{equation*}
  dZ_t^n = f(G_{\mu}^{-1}(t)) dt + \frac{\sigma}{\sqrt{n}}
  dW_t, \quad t \in [0,1],
\end{equation*}
where $W$ is a Brownian motion. Informally, if $\mu$ is known we
obtain by the time change $t = G_{\mu}(u)$,
\begin{equation*}
  d \wt Z_{u}^n = f(u) \mu(u) du + \sigma \sqrt{\frac{\mu(u)}{n}} d
  \wt W_u, \quad u \in [0,1],
\end{equation*}
where $\wt Z_u = Z_{G_{\mu}(u)}$ and $\wt W$ is a Brownian motion.
Finally, we obtain that~\eqref{eq:model_regression} is equivalent to
the heteroscedastic white noise model 
\begin{equation}
  \label{eq:hetero_white_noise_model}
  d Y_u^n = f(u) du + \frac{\sigma}{\sqrt{n \mu(u)}} dB_u, \quad u
  \in [0,1],
\end{equation}
where $B$ is a Brownian motion. In view of the result
by~\cite{donoho94} (see~\eqref{eq:donoho_result}) which is stated in
the model~\eqref{eq:white_noise_model} and comparing the noise levels
in the models~\eqref{eq:white_noise_model}
and~\eqref{eq:hetero_white_noise_model} (with $\sigma = 1$) we can
explain informally that our rate $r_{n, \mu}(\cdot)$ comes from the
former rate $\psi_n$ where we replace $n$ by $n \mu(x)$.

\subsection{About theorem \ref{thm:lower_bound}}
\label{sec:about_the_lower_bound}

From~\cite{bertin02}, we know when $s \in (0, 1]$ that
\begin{equation*}
  \liminf_n \inf_{\wh f_n} \sup_{f \in \Sigma} \Efm\big\{ w(
  v_{n, \mu}^{-1} \norminfty{\wh f_n - f}) \big\} \geq w(P),
\end{equation*}
where $v_{n, \mu}$ is given by~\eqref{eq:KB_rate_def}. An immediate
consequence is
\begin{equation}
  \label{eq:bad_lower_bound}
  \liminf_n\, \inf_{\wh f_n} \sup_{f \in \Sigma} \Efm\big\{ w \big(
  \sup_{x \in [0, 1]} r_{n, \mu}(x)^{-1} |\wh f_n(x) - f(x) | \big)
  \big\} \geq w(P),
\end{equation}
where it suffices to use $r_{n, \mu}(x) \leq v_{n, \mu}$ for any $x
\in [0, 1]$. This entails that $r_{n, \mu}(\cdot)$ is optimal in the
classical minimax sense, but this notion of optimality is weaker than
ours. Indeed, to prove the optimality of $r_{n, \mu}(\cdot)$ we need a
more "localised" version of the lower bound, hence
theorem~\ref{thm:lower_bound}.

In theorem \ref{thm:lower_bound}, if we choose $I_n = [0, 1]$ we find
back \eqref{eq:bad_lower_bound} and if $I_n = [\bar x - (\log
n)^{\gamma}, \bar x + (\log n)^{\gamma}] \cap [0, 1]$ for any $\gamma
> 0$ and $\bar x \in [0, 1]$ such that $\mu(\bar x) \neq \inf_{x \in
  [0, 1]} \mu(x)$, then obviously $v_{n, \mu}$ does not satisfy
\eqref{eq:stronger_lower_bound}.

\subsection{About proposition \ref{prop:CI}}

The confidence band $C_{n, \beta}(\cdot)$ is "design adaptive", in the
sense that it does not depend on $\mu$, but it depends on the
smoothness of $f$ via the parameters $s$ and $L$. The construction of
adaptive confidence bands is more involved. We know from~\cite{low97}
that the construction of an adaptive confidence band without extra
assumption is not feasible. However, if extra assumptions on the
smoothness of $f$ are supposed, it is possible to construct such
confidence bands,
see~\cite{picard_tribouley00},~\cite{hoffmann_lepski02}
and~\cite{cai_low04a, cai_low04b}. Here, we only focus on the
inhomogeneous aspect of the confidence band. Adaptation with respect
to the smoothness is beyond the scope of this study, and we would
encounter the same limitations.

\subsection{About assumption \ref{ass:design}}

In assumption~\ref{ass:design}, $\mu$ is supposed to be bounded from
below, and from above since it is continuous over $[0, 1]$. When $\mu$
is vanishing or exploding at a fixed point, we know from
\cite{gaiffas04a} that a wide range of pointwise minimax rates can be
achieved, depending on the behaviour of $\mu$ at this point. In this
case, we expect the optimal space dependent convergence rate (whenever
it exists) to be different from the classical minimax rate $\psi_n$
not only up to the constants but in order.

\section{Construction of an estimator}
\label{sec:estimator}

\subsection{Main idea}

The estimator $\wh f_n$ described below is using both kernel and local
polynomial methods. Its construction is divided in two parts: first,
at the discretisation points $x_j$ defined
by~\eqref{eq:discretisation_points}, we use a Nadaraya-Watson
estimator with a design data driven bandwidth. This part of the
estimator is used to attain the minimax constant. Between the
discretisation points, the estimator is defined by a Taylor expansion
where the derivatives estimates are done by local polynomial
estimation.

\subsection{The estimator at points $x_j$}

We consider the bandwidth $H_n(x)$ defined
by~\eqref{eq:data_driven_bandwidth_def} and we define
\begin{equation*}
  H_n^M = \max_{j \in \mc J_n} H_n(x_j),
\end{equation*}
where $x_j$ and $\mc J_n$ are defined in
section~\ref{sec:confidence_bands}. From~\cite{leonov97,leonov99} we
know that the function $\varphi_s$ defined
by~\eqref{eq:OR_optimal_function} is even and compactly supported.  We
denote by $[-T_s, T_s]$ its support and $\tau_n \eqdef \min(2 c_s T_s
H_n^M, \delta_n)$ where $\delta_n = (\log n)^{-1}$ and
\begin{equation}
  \label{eq:c_s_def}
  c_s \eqdef \Big( \frac{\sigma}{L} \Big)^{2 / (2s + 1)} \Big(
  \frac{2}{2s + 1}\Big)^{1 / (2s + 1)}.
\end{equation}
As usual with the estimation of a function over an interval, there is
a boundary correction. We decompose the unit interval into three parts
$[0, 1] = J_{n,1} \cup J_{n,2} \cup J_{n,3}$ where $J_{n,1} = [0,
\tau_n]$, $J_{n,2} = [\tau_n, 1 - \tau_n]$ and $J_{n,3} = [1 - \tau_n,
1]$. We also define $\mc J_{a,n} = \{ j | x_j \in J_{a,n} \}$ for $a
\in \{1, 2 ,3\}$. If $\varphi_s$ is defined
by~\eqref{eq:OR_optimal_function}, we consider the kernel
\begin{equation}
  \label{eq:OR_optimal_kernel}
  K_s = \frac{\varphi_s}{\int_{\setR} \varphi_s}.
\end{equation}
The "sharp" part of the estimator is defined as follows: at the points
$x_j$, we define $\wh f_n$ by
\begin{equation}
  \label{eq:estimator_def_xj}
  \wh f_n(x_j) \eqdef
  \begin{cases}
    \,\, \frac{ \displaystyle \frac{1}{n H_n(x_j)} \sumin Y_i K_s
      \Big( \frac{X_i - x_j}{c_s H_n(x_j)} \Big) }{\displaystyle \max
      \Big[ \delta_n, \, \frac{1}{n H_n(x_j) } \sumin K_s \Big(
      \frac{X_i - x_j}{c_s H_n(x_j)} \Big) \Big] } \,\, &\text{ if } j
    \in \mc J_{2,n}, \\ \,\, \bar{f}_n(x_j) \,\, &\text{ if } j \in
    \mc J_{1,n} \cup \mc J_{3,n}.
  \end{cases}
\end{equation}
This estimator is (up to the correction near the boundaries) a
Nadaraya-Watson estimator with the optimal kernel $K_s$ and a
bandwidth adjusted to the local amount of data. The boundary estimator
$\bar f_n$ is defined below.

\subsection{Between the points $x_j$ -- local polynomial estimation}
\label{sec:between_the_x_j}

We recall that $k = \ppint{s}$ where $s$ is the smoothness of the
unknown signal $f$. For any interval $I \subset [0, 1]$, we define the
inner product
\begin{equation*}
  \prodsca{f}{g}_I = \frac{1}{\bar \mu_n(I)} \int_I f g \, d\bar\mu_n,
\end{equation*}
where $\int_I f \, d \bar \mu_n = \sum_{X_i \in I} f(X_i) / n$. If $I
= I(x, h)$ -- see~\eqref{eq:I_x_h_def} -- for some $x \in [0, 1]$ and
$h > 0$, we define $\phi_{I, m}(y) = (y - x)^m$ and we introduce the
matrix $\mb X_I$ and vector $\mb Y_I$ with entries
\begin{equation*}
  (\mb X_I)_{p, q} = \prodsca{\phi_{I, p}}{\phi_{I, q}}_I \;
  \text{ and }\; (\mb Y_I)_{p} = \prodsca{Y}{\phi_{I, p}}_I,
\end{equation*}
for $0 \leq p, q \leq k$. Let us define
\begin{equation*}
  \bar{\mb X}_{I} = \mb X_{I} + \frac{1}{\sqrt{n \bar \mu_n(I)}} \, \mb
  I_{k+1} \, \ind{\Omega_{n, I}},
\end{equation*}
where $\Omega_{n, I} = \big\{ \lba(\mb X_{I}) \leq 1 / \sqrt{n \bar
  \mu_n(I)} \big\}$ and $\lba(M)$ is the smallest eigenvalue of a
matrix $M$ and $\mb I_{k+1}$ is the identity matrix on $\setR^{k+1}$.
Note that the correction term in $\bar{\mb X}_I$ entails
$\lba(\bar{\mb X}_I) \geq 1 / \sqrt{n \bar \mu_n(I)}$. When $\bar
\mu_n(I) > 0$, the solution $\wh \tta_I$ of the system
\begin{equation*}
  \bar{\mb X}_{I} \tta = \mb Y_{I},
\end{equation*}
is well defined. If $\bar \mu_n(I) = 0$, we take $\wh \tta_I = 0$.
Then, for any $1 \leq m \leq k$, a natural estimate of $f^{(m)}(x_j)$
is
\begin{equation*}
  \wt f_n^{(m)}(x_j) \eqdef m! (\wh \tta_{I(x_j, h_n)})_m,
\end{equation*}
where
\begin{equation*}
  h_n = (\sigma / L)^{2 / (2s + 1)} (\log n / n)^{1 / (2s + 1)},
\end{equation*}
and the estimator at the boundaries of $[0,1]$ is given by
\begin{equation*}
  \bar f_{n}(x_j) \eqdef (\wh \tta_{I(x_j, t_n)})_0,
\end{equation*}
where $t_n = (\sigma / L)^{2 / (2s + 1)} n^{-1 / (2s + 1)}$. Note that
the boundary estimator is a local polynomial estimator with the
pointwise bandwidth of estimation $t_n$. If we define
\begin{equation}
  \label{eq:Gamma_n_def}
  \Gamma_{n, I} = \Big\{ \min_{1 \leq m \leq k} \norm{\phi_{I, m}}_I
  \geq \frac{1}{\sqrt{n}} \Big\},
\end{equation}
where $\norm{\cdot}_I^2 = \prodsca{\cdot}{\cdot}_I$, then for $x \in
[x_j, x_{j+1})$, $j \in \mc J_n$, we take
\begin{equation}
  \label{eq:estimator_def_x}
  \wh f_n(x) \eqdef \wh f_n(x_j) + \Big( \sum_{m=1}^k \frac{\wt
    f_n^{(m)}(x_j)}{m!} (x - x_j)^{m} \Big) \ind{\Gamma_{n, I(x_j,
      h_n)}}.
\end{equation}

\section{Proof of theorem~1 and proposition~1}
\label{sec:proof_upper_bound}

The proof of theorem~\ref{thm:upper_bound} needs several preliminary
results. In section~\ref{sec:preparatory_results} we state the most
important lemmas while section~\ref{sec:LPA_results} is devoted to
useful results concerning local polynomial estimation. We delay the
proofs of these lemmas until section~\ref{sec:lemmas_proofs}, since
they can be skipped in a first reading. The proofs of
theorem~\ref{thm:upper_bound} and proposition~\ref{prop:CI} are given
in section~\ref{sec:proof_main_results_upper_bound}. We define the
risk
\begin{equation*}
  \mc E_{n, f} = \sup_{x \in [0, 1]} r_{n, \mu}(x)^{-1} |\wh f_n(x) -
  f(x)|,
\end{equation*}
and the discretised risk $\mc E_{n, f}^{\Delta} = \sup_{j \in \mc J_n}
r_{n, \mu}(x_j)^{-1} |\wh f_n(x_j) - f(x_j)|$. 

In the following, the notation $o(1)$ stands for a deterministic and
positive quantity going to $0$ as $n \raro +\infty$ indepedent of $f$
while $O(1)$ stands for a quantity bounded by a positive quantity
independent of $f$. If $A$ is non negative, we also define $O(A) =
O(1) \times A$. We denote $a \vee b = \max(a, b)$ and $a \wedge b =
\min(a, b)$. We consider the norms $\norminfty{g} = \sup_{x \in [0,
  1]} |g(x)|$, $\norm{g}_2 = (\int_0^1 g^2(x) dx)^{1/2}$, and
$\norminfty{x} = \max_{0 \leq m \leq k} |x_m|$, $\norm{x}_2 = (\sum_{0
  \leq m \leq k} x_m^2)^{1/2}$ when $x \in \setR^{k+1}$.

Since $\bar \mu_n(I(x, h)) / h$ is close to $\mu(x)$ in probability,
we have that $H_n(x)$ is close to
\begin{equation*}
  h_{n, \mu}(x) \eqdef \Big( \frac{\log n}{n \mu(x)} \Big)^{1 / (2s +
    1)}.
\end{equation*}
To avoid overloaded notations, it is convenient to write $K$ instead
of $K_s$ and to introduce for $j \in \mc J_n$,
\begin{equation*}
  H_j = H_n(x_j), \quad h_j = h_{n, \mu}(x_j), \quad \mu_j = \mu(x_j),
  \quad r_j = r_{n, \mu}(x_j),
\end{equation*}
\begin{equation*}
  K_{i, j} = K \Big( \frac{X_i - x_j}{c_s h_j} \Big), \quad \bar K_{i,
    j} = K \Big( \frac{X_i - x_j}{c_s H_j} \Big), \quad W_{i,j} =
  \frac{ \bar K_{i,j}}{\sumin \bar K_{i,j}},
\end{equation*}
and $q_j = n c_s h_j \mu_j$, $\bar q_j = n c_s H_j \mu_j$ where $c_s$
is given by~\eqref{eq:c_s_def}. We denote by $\mf X_n$ the sigma
algebra generated by the observations $X_i$, $1 \leq i \leq n$.

\subsection{Preparatory results}
\label{sec:preparatory_results}

We define
\begin{equation*}
  \mrm A_{n, j} \eqdef \Big\{ \big| \big( \sumin \bar K_{i,j} \big) /
  \bar q_j - 1 \big| \leq L_1 \delta_n^{s \wedge 1} \Big\},
\end{equation*}
where $L_1$ is a positive constant, and
\begin{align*}
  \mrm B_{n, j} &\eqdef \Big\{ \big| \big( \sumin K_{i,j} \big) / q_j
  - 1 \big| \leq \delta_n \Big\}, \quad \mrm C_{n, j} \eqdef \big\{
  | H_{j} / h_{j} - 1 | \leq \delta_n \big\}, \\
  \mrm E_{n, j} &\eqdef \Big\{ \big| \big( \sumin \bar K_{i,j}^2 \big)
  / q_j - \norm{K}_2^2 \big| \leq L_2 \delta_n^{s \wedge 1} \Big\},
\end{align*}
where $L_2$ is a fixed positive constant and
\begin{equation}
  \label{eq:mc_B_n_def}
  \mc B_n = \bigcap_{j \in \mc J_{2, n}} \big( \mrm A_{n, j} \cap \mrm
  B_{n, j} \cap \mrm E_{n, j} \big) \cap \bigcap_{j \in \mc J_n} \mrm
  C_{n, j}.
\end{equation}
A control over the probability of this event is given in
lemma~\ref{lem:awful_lemma} below.
Let us denote $Z_{n} = \max_{j \in \mc J_{2, n}} |Z_{n, j}|$ where
$Z_{n, j} = r_j^{-1} \sum_{i=1}^n \xi_i W_{i, j}$. Informally, the
variable $Z_{n}$ corresponds to the variance term of $\mc E_{n,
  f}^{\Delta}$. We recall that $M_n$ is equal to the cardinal of $\mc
J_n$.

\begin{lemma}[variance term]
  \label{lem:deviation_Znf}
  For any $\von > 0$,
  \begin{equation*}
    \sup_{f \in \Sigma^Q(s, L)} \Pfm \big \{ Z_{n} \ind{\mc B_n} > (1
    + \von) L c_s^s \norm{K}_2 \big \} \leq 2 (\log n)^{2s / (2s + 1)}
    n^{ -\von / (2s + 1)}.
  \end{equation*}
\end{lemma}

\begin{proof}
  Conditionally on $\mf X_n$, $Z_{n, j}$ is centered Gaussian with
  variance
  \begin{equation*}
    v_j^2 = \sigma^2 r_j^{-2} \sumin W_{ij}^2.
  \end{equation*}
  On $\mc B_n$, we have for any $j \in \mc J_{2, n}$ and $n$ large
  enough
  \begin{equation*}
    \sumin W_{i,j}^2 = \frac{ \sumin \bar K_{i,j}^2}{( \sumin \bar
      K_{i,j} )^2} \leq (1+o(1)) \frac{\norm{K}_2^2}{q_j}
    = (1 + o(1)) \frac{\norm{K}_2^2 r_j^2}{c_s \log n},
  \end{equation*}
  where we used the definition of $h_n(x)$, thus $v_j^2 \leq (1 +
  \von) \sigma^2 \norm{K}_2^2 / (c_s \log n)$. Using the standard
  Gaussian deviation, we obtain
  \begin{align*}
    \Pfm \{ |Z_{n, j}| \ind{\mc B_n} > (1 + \von) &L c_s^s \norm{K}_2
    \} \\
    &\leq 2 \exp \Big( -\frac{ (1 + \von) L^2 c_s^{2s + 1}}{2
      \sigma^2}
    \log n \Big) \\
    &= 2 \exp\Big( - \frac{(1 + \von)}{2s + 1} \log n \Big) = 2 n^{
      -(1 + \von) / (2s + 1) },
  \end{align*}
  and bounding from above the probability of $\cup_{j \in \mc J_{2,
      n}} \{ |Z_{n, j}| \ind{\mc B_n} > (1 + \von) L c_s^s \norm{K}_2
  \}$ by the sum of the probabilities, and since $|\mc J_{2, n}| \leq
  M_n \leq (\log n)^{2s / (2s + 1)} n^{1 / (2s + 1)}$, the lemma
  follows.
\end{proof}

For any $j \in \mc J_{n, 2}$, we define
\begin{equation*}
  b_{n, f} = \max_{j \in \mc J_{2,n}} | b_{n,f,j} | \quad \text{and}
  \quad U_{n, f} = \max_{j \in \mc J_{2,n}} | U_{n,f,j} |,
\end{equation*}
where $b_{n, f, j} = \Efm\{ B_{n,f,j} \ind{\mc B_n} \}$, $U_{n, f, j}
= B_{n,f,j} - b_{n,f,j}$ and
\begin{equation*}
  B_{n, f, j} = r_j^{-1} \sumin (f(X_i) - f(x_j)) W_{i,j}.  
\end{equation*}
The quantities $b_{n, f}$ and $U_{n, f}$ correspond to bias terms of
the risk $\mc E_{n, f}^{\Delta}$.

\begin{lemma}[first bias term]
  \label{lem:bias1}
  We have
  \begin{equation*}
    \limsup_n \sup_{f \in \Sigma(s, L) } b_{n, f} \leq L c_s^s \mc
    B(s, 1),
  \end{equation*}
  where $\mc B(s, L)$ is defined by \eqref{eq:B_s_L_def}.
\end{lemma}

\begin{lemma}[second bias term]
  \label{lem:bias2}
  There is a constant $D_U > 0$ such that for any $\von > 0$,
  \begin{equation*}
    \sup_{f \in \Sigma(s, L)} \Pfm \big\{ U_{n, f} \ind{\mc B_n} >
    \von \big \} \leq \exp \big( - D_U \, \von (1 \wedge \von) n^{2s /
      (2s + 1)} \big).
  \end{equation*} 
\end{lemma}

The proofs of these lemmas are delayed until section
\ref{sec:lemmas_proofs}.

\subsection{Local polynomial estimation}
\label{sec:LPA_results}

In this section we give results concerning local polynomial
estimation. This well known estimation procedure provides an efficient
method for recovering both a function and its derivatives. The
lemma~\ref{lem:biais_variance_derivatives_estimates} below is one
version of the bias variance decomposition of the local polynomial
estimator, which is classical: see~\cite{korostelev_tsybakov93},
\cite{fan_gijbels95,fan_gijbels96}, \cite{spok98} and
\cite{tsybakov03}, among many others. To a vector $\tta \in
\setR^{k+1}$ we associate the polynomial
\begin{equation*}
  P_{\tta}(y) = \tta_0 + \tta_1 y + \cdots + \tta_k y^k.
\end{equation*}
If $\wh \tta_I$ is the solution of the system $\bar{\mb X}_I \tta =
\mb Y_I$ (see section~\ref{sec:between_the_x_j}) for $I = I(x, h)$, we
define $\wh f_I(y) = P_{\wh \tta_I}(y - x)$. If $V_{I, k} = \Span\{
\phi_{I, m} ; 0 \leq m \leq k\}$, we note that on $\Omega_{n, I}$,
$\wh f_{I}$ satisfies
\begin{equation}
  \label{eq:LPA_variational}
  \prodsca{\wh f_{I}}{\phi}_{I} = \prodsca{Y}{\phi}_{I}, \quad \forall
  \phi \in V_{I,k}.
\end{equation}
By definition, we have $\wt f_n^{(m)}(x_j) = \wh f_{I(x_j,
  h_n)}^{(m)}(x_j)$, where $\wh f_{I}^{(m)}$ is the derivative of
order $m$ of $\wh f_{I}$, and $\bar f_{n}(x_j) = \wh f_{I(x_j,
  t_n)}(x_j)$, see section \ref{sec:between_the_x_j}. We introduce the
diagonal matrix $\mb \Lba_{I}$ with entries
\begin{equation*}
  (\mb \Lba_{I})_{m, m} = \norm{\phi_{I, m}}_{I}^{-1},
\end{equation*}
for $0 \leq m \leq k$, where $\norm{\cdot}_{I}^2 \eqdef
\prodsca{\cdot}{\cdot}_{I}$, the symmetrical matrix
\begin{equation*}
  \mc G_{I} \eqdef \mb \Lba_{I} \bar{\mb X}_{I} \mb \Lba_{I},
\end{equation*}
where $\bar{\mb X}_{I}$ is introduced in section
\ref{sec:between_the_x_j} and $\mc G$ the matrix with entries
\begin{equation*}
  (\mc G)_{p, q} = \frac{\chi_{p + q}}{\sqrt{\chi_{2p} \,\chi_{2q}}},
\end{equation*}
for $0 \leq p, q \leq k$, where $\chi_{m} = (1 + (-1)^{m}) / (2(m +
1))$. It is easy to see that $\lba(\mc G) > 0$ (we recall that
$\lba(M)$ is the smallest eigenvalue of a matrix $M$). We define the
event
\begin{equation*}
  \Omega_{n} = \bigcap_{j \in \mc J_n} \Omega_{n, I(x_j, h_n) }\cap
  \bigcap_{j \in \mc J_n} \Omega_{n, I(x_j, t_n)},
\end{equation*}
where $\Omega_{n, I}$ is defined in section \ref{sec:between_the_x_j}
and
\begin{equation*}
  \mc L_{n} = \bigcap_{j \in \mc J_n} \mc L_{n, I(x_j, h_n)} \cap
  \bigcap_{j \in \mc J_n} \mc L_{n, I(x_j, t_n)},
\end{equation*}
where if $I = I(x, h)$ for some $x \in [0, 1]$, $h > 0$,
\begin{equation*}
  \mc L_{n, I} = \{ | \lba(\mc G_I) - \lba(\mc G) | \leq \delta_n \}.
\end{equation*}
For $0 \leq m \leq 2k$ an interval $I \subset [0, 1]$ and $\delta >
0$, we define
\begin{equation*}
  \bar{\mrm D}_{n, m, I, \delta} \eqdef \bigg\{ \Big| \frac{1}{\bar
    \mu_n(I) |I|^m} \int_I \phi_{j, m} \, d \bar \mu_n - \chi_m \Big|
  \leq \delta \bigg\},
\end{equation*}
and
\begin{equation*}
  \mc D_n = \bigcap_{m=0}^{2k} \Big( \bigcap_{j \in \mc J_n} \bar{\mrm
    D}_{n, m, I(x_j, h_n), \delta_n} \cap \bigcap_{j \in \mc J_n}
  \bar{\mrm D}_{n, m, I(x_j, t_n), \delta_n} \Big).
\end{equation*}
%
We define
\begin{equation*}
  \mrm N_{n} = \bigcap_{j \in \mc J_n} \mrm N_{n, I(x_j, h_n)} \cap
  \bigcap_{j \in \mc J_n} \mrm N_{n, I(x_j, t_n)},
\end{equation*}
where
\begin{equation*}
  \mrm N_{n, I(x, h)} = \Big\{ \Big| \frac{\bar \mu_n(I(x, h))
  }{\mu(x) h} - 1 \Big| \leq \delta_n \Big\}.
\end{equation*}
Finally, we introduce
\begin{equation}
  \label{eq:mc_C_n_def}
  \mc C_n = \Omega_n \cap \mc L_n \cap \mc D_n \cap \mrm N_n.
\end{equation}
A control on the probability of this event is given in
lemma~\ref{lem:awful_lemma} below. We recall that $M_n$ is the
cardinal of $\mc J_n$.

\begin{lemma}
  \label{lem:biais_variance_derivatives_estimates}
  There exists a centered Gaussian vector $W \in \setR^{(k+1) M_n}$
  with
  \begin{equation*}
    \Efm\{ W_p^2 \}=1, \quad  0 \leq p \leq (k+1) M_n,
  \end{equation*}
  such that on $\mc C_n$, one has for any $0 \leq m \leq k$ and $f \in
  \Sigma(s, L)$\tup:
  \begin{equation}
    \label{eq:bias_variance_derivatives_estimates}
    \max_{j \in \mc J_n}| \wt f_n^{(m)}(x_j) - f^{(m)}(x_j) | \leq
    (1 + o(1)) C L h_n^{s - m} (1 + (\log n)^{-1/2} W^M),
  \end{equation}
  where 
  \begin{equation*}
    W^M \eqdef \max_{0 \leq p \leq (k+1) M_n} |W_p|,
  \end{equation*}
  and $C = C_{\lba, m, q, k}$ where $C_{\lba, m, q, k} = \lba^{-1}(\mc
  G) (k+1) m! \sqrt{2m+1} \big(1 \vee q^{-1/2} \big)$. For the
  estimator near the boundaries, we have for $a = 1$ and $a = 3$\tup:
  \begin{equation}
    \label{eq:bias_variance_boundary_estimator}
    \max_{j \in \mc J_{a,n}} |\bar f_n(x_j) - f(x_j) | \leq (1 + o(1))
    \bar C L t_n^{s} (1 + W^{(a)}),
  \end{equation}
  where 
  \begin{align*}
    W^{(1)} &= \max_{0 \leq p \leq (k+1) |\mc J_{1, n}|} |W_{p}| \\
    W^{(3)} &= \max_{(k+1)( |\mc J_{1,n}| + |\mc J_{2,n}|) + 1 \leq p
      \leq (k+1) M_n} |W_{p}|,
  \end{align*}
  and $\bar C = C_{\lba, 0, q, k}$.
\end{lemma}

\begin{lemma}
  \label{lem:LPA_bad_case}
  For any interval $I \subset [0, 1]$ and $p > 0$ we have
  \begin{equation*}
    \Efm\big\{ | (\wh \tta_I)_0 |^p  | \mf X_n \big\} = O(n^{p / 2}).
  \end{equation*}
  Moreover, for any $1 \leq m \leq k$, we have on $\Gamma_{n, I}$
  \tup(see section~\ref{sec:between_the_x_j}\tup)
  \begin{equation*}
    \Efm\big\{ | (\wh \tta_I)_m |^p  | \mf X_n \big\} = O(n^{p}).
  \end{equation*}
\end{lemma}

The proofs of these lemmas are delayed until section
section~\ref{sec:lemmas_proofs}. The following two lemmas are needed
for the proof of theorem~\ref{thm:upper_bound}.

\begin{lemma}
  \label{lem:risk_expectations_bounded}
  If $w(x) \leq A(1 + |x|)^b$ for some $A, b > 0$, we have
  \begin{equation}
    \label{eq:expectation_poynomially_bounded}
    \sup_{f \in \Sigma^Q(s, L)} \Efm\big\{ w^2(\mc E_{n, f} )
    \big\} = O\big(n^{2b (1 + s / (2s + 1)) }\big).
  \end{equation}
\end{lemma}

We define $\Gamma_n = \cap_{j \in \mc J_n} \Gamma_{n, I(x_j, h_n)}$
where $\Gamma_{n, I}$ is defined by~\eqref{eq:Gamma_n_def}. The
probability $\Pm$ stands for the joint law of the $X_1, \ldots, X_n$.

\begin{lemma}
  \label{lem:awful_lemma}
  There exists an event $\mc A_{n} \in \mf X_n$ such that for $n$
  large enough, under assumption~\ref{ass:design}
  \begin{equation}
    \label{eq:mc_A_n_deviation}
    \Pm\{ \mc A_{n}^c \} \leq \exp ( -D_{\mc A} n^{s / (2s + 1)}),
  \end{equation}
  where $D_{\mc A} > 0$ and
  \begin{equation}
    \label{eq:mc_A_n_subset_many_events}
    \mc A_{n} \subset \mc B_n \cap \mc C_n \cap \Gamma_n,
  \end{equation}
  where $\mc B_n$ is defined by~\eqref{eq:mc_B_n_def} and $\mc C_n$ is
  defined by~\eqref{eq:mc_C_n_def}.
\end{lemma}

\subsection{Proofs of the main results}
\label{sec:proof_main_results_upper_bound}

The next proposition is a deviation inequality for the discretised
risk $\mc E_{n, f}^{\Delta}$. This proposition is of special
importance in the proof of theorem~\ref{thm:upper_bound} and
proposition~\ref{prop:CI}.

\begin{proposition}
  \label{prop:main_risk_deviation}
  There is $D_{\mc E} > 0$ such that for any $\von > 0$, we have
  \begin{multline}
    \label{eq:deviation_discretized_risk}
    \sup_{f \in \Sigma^Q(s, L)} \Pfm\{ \mc E_{n, f}^{\Delta} \ind{\mc
      A_n} > (1 + \von) P \} \\
    \leq \exp\big( - D_{\mc E} \, \von( 1 \wedge \von) (\log n)^{2s /
      (2s + 1) } \big),
  \end{multline}
  for $n$ large enough. Moreover,
  \begin{equation}
    \label{eq:bounded_moment_discretized_risk}
    \sup_{f \in \Sigma^Q(s, L)} \Efm\big\{ w^2 (\mc E_{n, f}^{\Delta}
    \ind{\mc A_n}) \big\} = O(1).
  \end{equation}
\end{proposition}

\begin{proof}
  We decompose the risk into three parts
  \begin{equation}
    \label{eq:distretized_risk_decomposition}
    \mc E_{n, f}^{\Delta} = \mc E_{n, f}^{\Delta, 1} + \mc E_{n,
      f}^{\Delta, 2} + \mc E_{n, f}^{\Delta, 3},
  \end{equation}
  where $\mc E_{n, f}^{\Delta, a} = \sup_{j \in \mc J_{a, n}} r_j^{-1}
  | \wh f_n(x_j) - f(x_j)|$. For $a=1$ and $a= 3$, the quantity $\mc
  E_{n, f}^{\Delta, a}$ is the risk at the boundaries of $[0, 1]$.
  Note that on $\mc B_n$, we have $\sumin \bar K_{i, j} / (n H_j) >
  c_s \mu_j (1 - L_1 \delta_n^{s \wedge 1}) > c_s q( 1 - L_1
  \delta_n^{s \wedge 1}) > \delta_n$ for $n$ large enough. Hence,
  since $\mc A_n \subset \mc B_n$ (see lemma~\ref{lem:awful_lemma}) we
  can decompose on $\mc A_n$ the middle risk into bias and variance
  terms as follows:
  \begin{equation}
    \label{eq:Efn2_decomp}
    \mc E_{n, f}^{\Delta, 2} \leq b_{n, f} + U_{n, f} + Z_{n}.
  \end{equation}
  In view of lemma \ref{lem:bias1} we have for $n$ large enough $b_{n,
    f} \leq (1 + 2\von) L c_s^s \mc B(s, 1)$ and using equation
  \eqref{eq:P_s_decomposition} we obtain
  \begin{align*}
    \{ \mc E_{n, f}^{\Delta, 2} \ind{\mc A_n} &> (1 + 2 \von) P \} \\
    &\subset \{ Z_n \ind{\mc B_n} > (1 + \von) L c_s^s \norm{K}_2 \}
    \cup \{ U_{n, f} \ind{\mc B_n} > \von L c_s^s \norm{K}_2 \}.
  \end{align*}
  Then, in view of the lemmas~\ref{lem:deviation_Znf} and
  \ref{lem:bias2}, it is easy to find $D_{2} > 0$ such that for any $f
  \in \Sigma^Q(s, L)$ and $n$ large enough,
  \begin{equation}
    \label{eq:E_n_f_2_deviation}
    \Pfm\big\{ \mc E_{n,f}^{\Delta, 2} \ind{\mc A_n} > (1 + 2 \von)
    P \big\} \leq \exp\big( - D_2 \, \von(1 \wedge \von) \log n
    \big).
  \end{equation}
  Using lemma~\ref{lem:biais_variance_derivatives_estimates}, we
  obtain
  \begin{equation}
    \label{eq:E_nf_1_majoration}
    \mc E_{n, f}^{\Delta, 1} \ind{\mc A_n} \leq L_3 \delta_{n}^{s / (2
      s + 1)}(1 + W^{(1)}),
  \end{equation}
  where $W^{(1)} = \max_{0 \leq p \leq (k+1) \times |\mc J_{1, n}|}
  |W_p|$ and $L_3 = \bar C \norminfty{\mu}^{s / (2s + 1)}$. Since $W$
  is a centered Gaussian vector such that $\Efm\{ W_p^2 \} = 1$ for $0
  \leq p \leq (k+1) M_n$ it is well known (see for instance
  in~\cite{ledoux_talagrand91}) that
  \begin{equation*}
    \Efm\{ W^{(1)} \} \leq \sqrt{2 \log( (k+1) |\mc J_{n, 1}|) } = O(
    \sqrt{\log \log n}),
  \end{equation*}
  since $|\mc J_{1,n}| = O(\log n)$, and that for any $\lba > 0$,
  \begin{equation*}
    \Pfm\big\{ W^{(1)} - \Efm\{ W^{(1)}  \}  > \lba \big\} \leq 2\exp(
    - \lba^2 / 2).
  \end{equation*}
  Then, when $n$ is large enough,
  \begin{align*}
    \Pfm\big\{ \mc E_{n, f}^{\Delta, 1} \ind{\mc A_n} > 2\von P \big\}
    &\leq \Pfm\big\{ W^{(1)} - \Efm\{ W^{(1)} \} > \von
    P \delta_n^{-s / (2s+1)} / L_3 \big\} \\
    &\leq 2 \exp\big( - \von^2 P^2 \delta_n^{ -2s / (2s + 1)} / (2
    L_3^2) \big).
  \end{align*}
  The same result holds for $\mc E_{n, f}^{\Delta, 3}$. Hence,
  together with~\eqref{eq:E_n_f_2_deviation}, for a good choice of
  $D_{\mc E}$ we obtain~\eqref{eq:deviation_discretized_risk}. It is
  easy to prove~\eqref{eq:bounded_moment_discretized_risk}
  from~\eqref{eq:deviation_discretized_risk}. For any $f \in
  \Sigma^Q(s, L)$ and $p > 0$, when $n$ is large enough,
  \begin{align*}
    \Efm\{ (\mc E_{n, f}^{\Delta})^p \ind{\mc A_n} \} &= p
    \int_{0}^{+\infty} t^{p - 1} \Pfm\{ \mc E_{n, f}^{\Delta} \ind{\mc
      A_n} > t \} dt \\
    &\leq (2 P)^p + p e^{D_{\mc E}} \int_{2 P}^{+\infty} t^{p - 1}
    \exp\big( -D_{\mc E} t / P \big) dt = O(1),
  \end{align*}
  thus~\eqref{eq:bounded_moment_discretized_risk}, since $w(x) \leq
  A(1 + |x|^b)$.
\end{proof}

\begin{proof}[Proof of theorem \ref{thm:upper_bound}]
  Let $x \in [x_j, x_{j+1})$. Since $\mu \in \Sigma(\nu, \varrho)$
  with $0 < \nu \leq 1$ we have clearly $\mu^{s/(2s+1)} \in
  \Sigma(s\nu / (2s+1), \varrho^{s/(2s+1)})$ and using assumption
  \ref{ass:design},
  \begin{equation}
    \label{eq:rnx_minus_rnxk}
    \sup_{x \in [x_j, x_{j+1}]} | r_{n, \mu}(x)^{-1} - r_j^{-1} | \leq
    r_j^{-1} \Big( \frac{\varrho}{q} \Big)^{s / (2s + 1)} \Delta_n^{s
      \nu / (2s + 1)} = o(1) r_j^{-1}.
  \end{equation}
  Since $f \in \Sigma^Q(s, L)$, writing the Taylor expansion of $f$ at
  $x \in [x_j, x_{j+1})$ we obtain:
  \begin{align*}
    |\wh f_n(x) - f(x)| &\leq |\wh f_n(x_j) - f(x_j)| \\
    &+ \sum_{m=1}^k ( \wt f_n^{(m)}(x_j) - f^{(m)}(x_j)) \frac{(x -
      x_j)^m}{m!} + L \Delta_n^s,
  \end{align*}
  and in view of~\eqref{eq:rnx_minus_rnxk},
  \begin{equation*}
    \mc E_{n, f} \leq (1 + o(1)) \Big( \mc E_{n, f}^{\Delta} + \max_{j
      \in \mc J_n} r_j^{-1} \sum_{m=1}^k | \wt f_n^{(m)}(x_j) -
    f^{(m)}(x_j) | \frac{\Delta_n^m}{m!} \Big) + O(\delta_n^s).
  \end{equation*}
  We consider the event $\mc A_n$ from lemma~\ref{lem:awful_lemma}.
  Since $\mc A_n \subset \mc C_n$ we have that on $\mc A_n$, in view
  of lemma~\ref{lem:biais_variance_derivatives_estimates} and for any
  $1 \leq m \leq k$,
  \begin{align*}
    \max_{j \in \mc J_n} r_j^{-1} |\wt f_n^{(m)}(x_j) &- f^{(m)}(x_j)|
    \frac{\Delta_n^m}{m!} \\
    &\leq (1 + o(1)) \delta_n^m \norminfty{\mu}^{s / (2s+1)} C (1 +
    (\log n)^{-1/2} W^M),
  \end{align*}
  and then
  \begin{equation*}
    \mc E_{n, f} \ind{\mc A_n} \leq (1 + o(1)) \mc E_{n, f}^{\Delta}
    \ind{\mc A_n} + O(1) \delta_n (1 + \delta_n^{1/2} W^M) + o(1).
  \end{equation*}
  We define $\mc W_n \eqdef \{ |W^M - \Efm\{ W^M \}| \leq
  \delta_n^{-1} \}$. Since $W^M = \max_{0 \leq p \leq (k+1) M_n}
  |W_p|$, we know in the same way as in the proof of
  proposition~\ref{prop:main_risk_deviation} that $\Efm\{ W^M \} \leq
  \sqrt{2 \log((k+1) M_n)} = O(\delta_n^{-1/2})$ and
  \begin{equation}
    \label{eq:W_n_comp_dev}
    \Pfm \{ \mc W_n^c \} \leq 2 \exp( - \delta_n^{-2} / 2).
  \end{equation}
  Thus
  \begin{equation}
    \label{eq:risk_bounded_by_discretized_risk}
    \mc E_{n, f} \ind{\mc A_n \cap \mc W_n} \leq (1+o(1)) \mc E_{n,
      f}^{\Delta} \ind{\mc A_n} + o(1),
  \end{equation}
  and since $w$ is non-decreasing, we have for any $\von > 0$
  \begin{align*}
    \Efm\{ w&(\mc E_{n, f}) \} \\
    &\leq \Efm\{ w(\mc E_{n, f}) \ind{\mc A_n \cap \mc W_n}\} + \Efm\{
    w(\mc E_{n, f}) \ind{\mc A_n^c \cup \mc W_n^c}\} \\
    &\leq w ((1 + 2 \von) P) + \big( \Efm\{ w^2(\mc E_{n, f}) \} \,
    \Pfm\{ \mc A_n^c \cup \mc W_n^c \} \big)^{1/2}  \\
    &+ \big( \Efm \big\{ w^2 \big( (1 + 2\von) \mc E_{n, f}^{\Delta}
    \ind{\mc A_n} \big) \big\} \, \Pfm\{ \mc E_{n, f}^{\Delta}
    \ind{\mc A_n} > (1 + \von)P \} \big)^{1/2} \\
    &\leq w ((1 + 2 \von) P) + O \big(n^{b (1 + s / (2s + 1))}
    \exp(-(\log n)^2 / 4) \big) \\
    &+ O \big( \exp(- D_{\mc E} \, \von(1 \wedge \von) (\log n)^{2s /
      (2s + 1)}) \big) = w ((1 + 2 \von) P) + o(1),
  \end{align*}
  where we used proposition~\ref{prop:main_risk_deviation},
  lemmas~\ref{lem:risk_expectations_bounded},~\ref{lem:awful_lemma}
  and the fact that $w$ is continuous. Thus,
  \begin{equation*}
    \limsup_n \sup_{f \in \Sigma^Q(s, L)} \Efm\{ w(\mc E_{n, f}) \}
    \leq w((1 + 2\von) P),
  \end{equation*}
  which concludes the proof of theorem~\ref{thm:upper_bound} since
  $\von$ can be chosen arbitrarily small.
\end{proof}

\begin{proof}[Proof of proposition~\ref{prop:CI}]
  We consider the event $\mc W_n$ defined in the proof of
  theorem~\ref{thm:upper_bound}. Since $\mc A_n \subset \mc B_n
  \subset \mrm C_{n, j}$ for any $j \in \mc J_n$ we have
  \begin{equation}
    \label{eq:R_j_and_r_j}
    (1 - o(1)) r_j \leq R_n(x_j) \leq (1 + o(1)) r_j
  \end{equation}
  on $\mc A_n$. In view of~\eqref{eq:rnx_minus_rnxk}
  and~\eqref{eq:risk_bounded_by_discretized_risk} we have for any $j
  \in \mc J_n$, $x \in [x_j, x_{j+1})$ on $\mc A_n \cap \mc W_n$
  \begin{align*}
    R_n(x)^{-1} | \wh f_n(x) - f(x) | &= \frac{r_{n,
        \mu}(x)}{R_n(x_j)} r_{n, \mu}(x)^{-1} |\wh f_n(x) - f(x)| \\
    &\leq (1 + o(1)) \mc E_{n, f} \leq (1 + o(1)) \mc E_{n,
      f}^{\Delta} + o(1).
  \end{align*}
  Thus, if $\mc F_{n, f, \beta} = \big\{ \sup_{x \in [0, 1]}
  R_n(x)^{-1} |\wh f_n(x) - f(x) | \leq (1 + \beta) P \big\}$
  lemma~\ref{lem:awful_lemma},
  proposition~\ref{prop:main_risk_deviation}
  and~\eqref{eq:W_n_comp_dev} entail for any $f \in \Sigma^Q(s, L)$,
  \begin{align*}
    \Pfm\{ \mc F_{n, f, \beta}^c \} &\leq \Pfm\{ \mc F_{n, f, \beta}^c
    \cap \mc A_n \cap \mc W_n \} + \Pfm\{ \mc A_n^c \cup \mc W_n^c \}
    \\
    &\leq \Pfm\{ \mc E_{n, f}^{\Delta} \ind{\mc A_n} > (1 + \beta / 2)
    P \} + \Pfm\{ \mc A_n^c \cup \mc W_n^c \} \\
    &\leq \exp( -D_{c} \, \beta(2 \wedge \beta) (\log n)^{2s / (2s +
      1)} \big),
  \end{align*}
  for a good choice of $D_c$. When $n$ is large enough, the choice
  $\beta = \beta(n, \alpha)$ makes the last part of the above
  inequality equal to $\alpha$, hence~\eqref{eq:CI_coverage_prob}.
  Using again~\eqref{eq:R_j_and_r_j}, lemma~\ref{lem:awful_lemma}
  and~\eqref{eq:rnx_minus_rnxk} it is easy to
  obtain~\eqref{eq:CI_length_equivalent}.
\end{proof}

\subsection{Proof of lemmas~\ref{lem:bias1}, \ref{lem:bias2},
  \ref{lem:biais_variance_derivatives_estimates},
  \ref{lem:LPA_bad_case}, \ref{lem:risk_expectations_bounded} and
  \ref{lem:awful_lemma}}
\label{sec:lemmas_proofs}

Since $b_{n, f}$ and $U_{n, f}$ only depend on $f$ via its values in
$[0, 1]$, we have
\begin{equation}
  \label{eq:sup_over_01_sup_over_R}
  \sup_{f \in \Sigma(s, L)} b_{n, f} = \sup_{f \in \Sigma(s, L;
    \setR)} b_{n, f}, \quad \sup_{f \in \Sigma(s, L)} U_{n, f} =
  \sup_{f \in \Sigma(s, L; \setR)} U_{n, f}.
\end{equation}
Here, it is convenient to introduce $P_{j} \eqdef \sumin ( f(X_i) -
f(x_j) ) \bar K_{i, j}$ and $Q_{j} \eqdef \sumin \bar K_{i, j}$.

\begin{proof}[Proof of lemma~\ref{lem:bias1}]
  On $\mrm A_{n, j} \cap \mrm C_{n, j}$ we have $(1 - o(1)) q_j \leq
  Q_j \leq (1 + o(1)) q_j$ and since $\mc B_{n} \subset \mrm A_{n, j}
  \cap \mrm C_{n, j}$ for any $j \in \mc J_{2, n}$, we have
  \begin{equation*}
    | b_{n, f, j} | = r_j^{-1} | \Efm \{ (P_j / Q_j) \ind{\mc B_n} \}
    | \leq  (1+o(1)) (r_j q_j)^{-1} | \Efm\{ P_{j} \ind{\mc B_n} \} |.
  \end{equation*}
  Recalling that $K = \varphi_s / \int \varphi_s$ with $\varphi_s \in
  \Sigma(s, 1; \setR)$ we have for any $x, y \in \setR$
  \begin{equation*}
    |K(x) - K(y)| \leq \kpa |x - y|^{s_1},
  \end{equation*}
  where $s_1 = s \wedge 1$ and $\kpa = (\int \varphi_s)^{-1}$ when $s
  \in (0, 1]$ and $\kpa = \norminfty{K'}$ when $s > 1$. Since $\supp K
  = [-T_s, T_s]$, we have for $n$ large enough on $\mc B_n$:
  \begin{equation}
    \label{eq:use_of_K_regularity}
    \begin{split}
      | \bar K_{i, j} - K_{i, j} | &\leq \kpa \Big| \frac{X_i -
        x_j}{c_s H_j} \Big|^{s_1} \Big|\frac{H_j}{h_j} - 1 \Big|^{s_1}
      \ind{|X_i - x_j| \leq c_s T_s (H_j \vee h_j)} \\
      &\leq \kpa T_s^{s_1} \Big( \frac{\delta_n}{1 - \delta_n}
      \Big)^{s_1} \ind{|X_i - x_j| \leq c_s T_s (1+\delta_n) h_j} =
      o(1) \ind{\mrm M_{i,j}},
    \end{split}
  \end{equation}
  where $\mrm M_{i, j} \eqdef \{|X_i - x_j| \leq c_s T_s (1+\delta_n)
  h_j \}$. We introduce $\nu_{f,j}(x) = \ind{f(x) \geq f(x_j)} -
  \ind{f(x) < f(x_j)}$, $R_{i,j} = (f(X_i) - f(x_j)) K_{i,j}$,
  $S_{i,j} = \nu_{f, j}(X_i) (f(X_i) - f(x_j)) \ind{\mrm M_{i,j}}$,
  $R_j = \sumin R_{i,j}$ and $S_j = \sumin S_{i,j}$. Then,
  \begin{align*}
    &\frac{1}{r_j q_j} | \Efm\{ P_{j} \ind{\mc B_n }\} \big|  \\
    &\leq \frac{1}{r_j q_j} \big( | \Efm\{ R_j \} | + o(1) | \Efm\{
    S_j \} | \big) \\
    &\leq \frac{1}{r_j \mu_j} \Big( \big| \int (f(x_j + y c_s h_j) -
    f(x_j)) K(y) \mu(x_j + y c_s h_j) dy \big| \\
    &+ o(1) \big| \int_{|y| \leq (1 + \delta_n) T_s} (f(x_j + y c_s
    h_j) - f(x_j)) \nu_{f, j}(x_j + c_s y h_j) \mu(x_j + y c_s h_j) dy
    \big| \Big),
  \end{align*}
  and since $\mu \in \Sigma_q(\nu, \varrho)$ we have
  \begin{align*}
    b_{n, f, j} &\leq \frac{1 + o(1)}{r_j} \big| \int (f(x_j + y c_s
    h_j) - f(x_j)) K(y) dy \big| \\
    &+ \frac{o(1)}{r_j q} \int_{|y| \leq 2 T_s} |f(x_j + y c_s h_j) -
    f(x_j)| dy.
  \end{align*}
  Using~\eqref{eq:sup_over_01_sup_over_R} and the fact that $\Sigma(s,
  L; \setR)$ is invariant by translation,
  \begin{multline}
    \sup_{f \in \Sigma(s, L; \setR)} b_{n,f, j} \leq (1 + o(1))
    \sup_{f \in \Sigma(s, L; \setR)} \max_{j \in \mc J_{2, n}}
    \frac{1}{r_j} \Big( \big| \int (f(c_s h_j y) - f(0)) K(y) dy \big| \\
    + o(1) \int_{|y| \leq 2T} |f(c_s h_j y) - f(0)| dy \Big).
  \end{multline}
  Now we use an argument which is known as \emph{renormalisation}, see
  \cite{donoho_low_92}. We introduce the functional operator $\mc
  U_{a, b}f(\cdot) = af(b\,\cdot)$. We have that $f \in \Sigma(s, L;
  \setR)$ is equivalent to $\mc U_{a, b} f \in \Sigma(s, L a b^s;
  \setR)$. Then, choosing $a = (L c_s^s h_j^s)^{-1}$ and $b = c_s h_j$
  entails
  \begin{equation*}
    \sup_{f \in \Sigma(s, L; \setR)} b_{n, f} \leq (1 + o(1)) L c_s^s
    \mc B(s, 1) + o(1) \sup_{f \in \Sigma(s, 1; \setR)} \int_{|y| \leq
      2T} |f(y) - f(0)| dy,
  \end{equation*}
  where $\mc B(s, 1)$ is given by~\eqref{eq:B_s_L_def} and where we
  recall that $r_j = h_j^s$. We define $f_{k}(y) = f(0) + f^{'}(0) y +
  \cdots + f^{(k)}(0) y^k / k!$. Since $f \in \Sigma(s, L; \setR)$, we
  have $f - f_{k} \in \Sigma(s, L; \setR)$ and finally
  \begin{equation*}
    \sup_{f \in \Sigma(s, L; \setR)} b_{n, f} \leq (1 + o(1)) L c_s^s
    \mc B(s, 1) + o(1) \int_{|y| \leq 2T} |y|^{s} dy . \qedhere
  \end{equation*}
\end{proof}

\begin{proof}[Proof of lemma~\ref{lem:bias2}]
  We recall that $U_{n, f, j} \eqdef r_j^{-1} ( B_j - \Efm\{ B_j
  \ind{\mc B_n } \} )$. We use the same notations as in the proof of
  lemma~\ref{lem:bias1}. On $\mc B_n$ we have $(1 - o(1)) q_j \leq Q_j
  \leq (1 + o(1)) q_j$, and since $\Efm\{ P_j^2\} \leq 4 Q^2
  \norminfty{K}^2 n^2$ we obtain in view of
  lemma~\ref{lem:awful_lemma}:
  \begin{equation*}
    \frac{1}{r_j q_j} | \Efm\{ P_j \ind{\mc B_n^c} \} | \leq
    \frac{1}{r_j q_j} \sqrt{\Efm\{ P_j^2 \}} \sqrt{\Pm\{ \mc B_n^c \}}
    = o(1).
  \end{equation*}
  Then, it is easy to see that on $\mc B_n$,
  \begin{equation*}
    |U_{n, f, j}| \leq \frac{1}{r_j q_j} \Big( (1 + o(1)) \big| P_j -
    \Efm\{ P_j \} \big| + o(1) \big| \Efm\{ P_j \ind{\mc B_n} \} \big|
    \Big) + o(1),
  \end{equation*}
  and we know from the proof of lemma~\ref{lem:bias1} that
  \begin{equation*}
    \frac{1}{r_j q_j} | \Efm\{ P_j \ind{\mc B_n} \} | \leq \sup_{f \in
      \Sigma(s, L)} \max_{j \in \mc J_{2, n}} \frac{1}{r_j
      q_j} | \Efm\{ P_j \ind{\mc B_n} \} | \leq (1 + o(1)) L c_s^s \mc
    B(s, 1),
  \end{equation*}
  thus $| U_{n, f, j} | \leq (1 + o(1)) (r_j q_j)^{-1} | P_j - \Efm\{
  P_j \} | + o(1)$ on $\mc B_n$. From the proof of lemma
  \ref{lem:bias1}, we know that $(r_j q_j)^{-1} |\Efm\{ S_j \}| =
  O(1)$, and using~\eqref{eq:use_of_K_regularity} it is an easy
  computation to obtain that on $\mc B_n$,
  \begin{equation*}
    |P_j - \Efm\{ P_j \}| \leq |R_j - \Efm\{ R_j \}| + o(1) |S_j -
    \Efm\{ S_j \}| + o(1) |\Efm\{ S_j \}|.
  \end{equation*}
  Then we have for $n$ large enough
  \begin{multline*}
    \Pfm\{ |U_{n, f, j}| \ind{\mc B_n} > \von \} \leq \Pfm\big\{ | R_j
    - \Efm\{ R_j \} | > \frac{\von r_j q_j}{3} \big\} \\ + \Pfm\big\{
    | S_j - \Efm\{ S_j \} | > \frac{\von r_j q_j}{3} \big\}.
  \end{multline*}
  We use Bernstein inequality to the sum of variables $\bar R_{i,j}
  \eqdef R_{i,j} - \Efm\{ R_{i,j} \}$ and $\bar S_{i,j} \eqdef S_{i,j}
  - \Efm\{ S_{i,j} \}$, $1 \leq i \leq n$. The variables $(\bar
  R_{i,j})_{1 \leq i \leq n}$ are clearly independent, centered and
  satisfy $| \bar R_{i,j} | \leq 4 Q K_{\infty}$. In view
  of~\eqref{eq:sup_over_01_sup_over_R} and since $\mu \in
  \Sigma_q(\nu, \varrho)$, it is easy to prove with the same arguments
  as in the end of the proof of lemma~\ref{lem:bias1} that
  \begin{align*}
    \Efm\{ \bar R_{i,j}^2 \} &\leq \Efm\{ R_{i,j}^2 \} \\
    &\leq (1 + o(1)) c_s h_j \mu_j \int (f(x_j + c_s h_j y) -
    f(x_j))^2 K^2(y) dy \\
    &\leq (1 + o(1)) c_s h_j \mu_j \sup_{f \in \Sigma(s, L; \setR)}
    \int (f(x_j + c_s h_j y) - f(x_j))^2 K^2(y) dy \\
    &\leq (1 + o(1)) L^2 (c_s h_j)^{2s+1} \mu_j \sup_{f \in \Sigma(s,
      L; \setR)} \int (f(y) - f(0))^2 K^2(y) dy \\
    &\leq (1 + o(1)) L^2 (c_s h_j)^{2s+1} \mu_j \int y^{2s} K^2(y) dy
    / (k!)^2.
  \end{align*}
  Then $\sumin \Efm\{ \bar R_{i,j}^2 \} = O( r_j^2 q_j )$ and the
  Bernstein inequality entails that for $n$ large enough, there is a
  constant $D_4 > 0$ such that
  \begin{align*}
    \Pfm\{ |R_j - \Efm\{ R_j \}| > \von r_j q_j / 3 \} \leq 2 \exp(
    -D_4 \von (1 \wedge \von) n^{s / (2s + 1) } ).
  \end{align*}
  The variables $(\bar S_{i,j})_{1 \leq i \leq n}$ are independent,
  centered and such that $|\bar S_{i,j}| \leq 4 Q$, and in the same
  way as previously we can prove $\sumin \Efm\{ \bar S_{i,j}^2 \} =
  O(r_j^2 q_j)$. Using again Bernstein inequality, it is easy to find
  $D_5$ such that
  \begin{align*}
    \Pfm\{ |S_j - \Efm\{ S_j \}| > \von r_j q_j / 3 \} \leq 2 \exp(
    -D_5 \von( 1 \wedge \von) n^{s / (2s + 1)} ),
  \end{align*}
  and since $|\mc J_{2, n}| \leq M_n$, we have for any $f \in
  \Sigma^Q(s, L)$,
  \begin{align*}
    \Pfm\{ |U_{n, f}| \ind{\mc B_n} > \von \} &\leq \sum_{j \in \mc
      J_{2,n}} \Pfm\{ |U_{n, f, j}| \ind{\mc B_n} > \von \} \\
    &\leq 4 M_n \exp \big( - (D_4 \wedge D_5) \, \von(1 \wedge \von)
    n^{s / (2s + 1)} \big).
  \end{align*}
  Since $4 M_n \exp( -(D_4 \wedge D_5) \von (1 \wedge \von) n^{s / (2s
    + 1)} / 2)$ goes to $0$ as $n$ goes to $+\infty$, the lemma
  follows with $D_U = (D_4 \wedge D_5) / 2$.
\end{proof}

\begin{proof}[Proof of
  lemma~\ref{lem:biais_variance_derivatives_estimates}]
  \label{proof:lema_bias_variance_chap3}
  We take $I = I(x, h)$ for some $x \in [0, 1]$, $h > 0$ and define
  the vector $\tta_I$ with coordinates $(\tta_I)_m = f^{(m)}(x) / m!$
  for $0 \leq m \leq k$.  Since $\bar{\mb X}_{I} = \mb X_{I}$ on
  $\Omega_{n, I}$, we have $\mb \Lba_{I}^{-1} (\wh \tta_{I} - \tta_I)
  = \mc G_{I}^{-1} \mb \Lba_{I} \mb X_{I} (\wh \tta_{I} - \tta_I)$. If
  $f_{I}(y) = P_{\tta_I}(y - x)$, we have in view
  of~\eqref{eq:LPA_variational} for any $0 \leq m \leq k$:
  \begin{align*}
    (\mb X_{I}(\wh \tta_{I} - \tta_I))_m &= \prodsca{\wh f_{I} -
      f_{I}}{\phi_{I, m}}_{I} = \prodsca{Y - f_{I}}{\phi_{I, m}}_{I} \\
    &= \prodsca{f - f_{I}}{\phi_{I, m}}_{I} + \prodsca{\xi}{\phi_{I,
        m}}_{I},
  \end{align*}
  thus $\mb X_{I} (\wh \tta_{I} - \tta_I) \eqdef \mb B_{I} + \mb
  V_{I}$.  Since $f \in \Sigma(s, L)$,
  \begin{equation*}
    (\mb \Lba_{I} \mb B_{I})_m \leq \norm{\phi_{I, m}}_{I}^{-1} |
    \prodsca{f - f_{I}}{\phi_{I, m}}_{I} | \leq \norm{f - f_{I}}_{I}
    \leq L h^s / k!,
  \end{equation*}
  then we can write
  \begin{equation*}
    \mb \Lba_{I}^{-1} (\wh \tta_{I} - \tta_I) = \mc G_{I}^{-1} \frac{L
      h^s}{k!} u + \frac{\sigma}{\sqrt{n \bar \mu_n(I)}} \, \mc
    G_{I}^{-1/2} \gamma_{I},
  \end{equation*}
  where $u \in \setR^{k+1}$ is such that $\norminfty{u} \leq 1$ and
  $\gamma_{I} = (\sigma \sqrt{n \bar \mu_n(I)})^{-1} \mc G_{I}^{-1/2}
  \mb \Lba_{I} \mb D_{I} \xi \eqdef \mb T_{I} \xi$, where $\mb D_{I}$
  is the matrix of size $n \bar \mu_n(I) \times (k+1)$ with entries
  $(\mb D_{I})_{i,m} = (X_i - x)^m$, so that $\mb X_{I} = (n \bar
  \mu_n(I))^{-1} \mb D_{I}' \mb D_{I} $. Since $\mb T_I' \mb T_I =
  \sigma^{-1} \mb I_{k+1}$, we obtain that $\gamma_{I}$ is,
  conditionally on $\mf X_n$, centered Gaussian with covariance equal
  to $\mb I_{k+1}$.

  Consider $I = I(x_j, h)$ for some $j \in \mc J_n$, $h > 0$. From the
  inequality $\norminfty{\cdot} \leq \norm{\cdot} \leq \sqrt{k+1}
  \norminfty{\cdot}$ and since $\norm{\mc G_{I}^{-1/2}} \leq
  \sqrt{k+1} \norm{\mc G_{I}^{-1}}$ ($\mc G_{I}$ is symmetrical with
  entries smaller than $1$ in absolute value) we get
  \begin{align*}
    \norminfty{ \mb \Lba_{I}^{-1} (\wh \tta_{I} - \tta_I) } &\leq
    \norminfty{\mc G_{I}^{-1} \frac{L h^s}{k!} u} +
    \frac{\sigma}{\sqrt{n \bar \mu_n(I)}} \norminfty{\mc G_{I}^{-1/2}
      \gamma_{I}} \\
    &\leq \norm{\mc G_{I}^{-1}} (k+1) \big( L h^s +
    \frac{\sigma}{\sqrt{n \bar \mu_n(I) }} \norminfty{\gamma_{I}}
    \big) \\
    &= \lba^{-1}(\mc G_{I}) (k+1) \big( L h^s + \frac{\sigma}{\sqrt{n
        \bar \mu_n(I)}} \max_{0 \leq m \leq k} |W_{(k+1) j + m}|
    \big),
  \end{align*}
  where $W \eqdef (\gamma_{I(x_0, h)}, \ldots, \gamma_{I(x_{M_n}, h)
  })'$. If $\mb T \eqdef (\mb T_{I(x_0, h)}, \ldots, \mb T_{I(x_{M_n},
    h)})'$ we have $W = \mb T \xi$, thus $W$ is a centered Gaussian
  vector and for any $(k+1)j \leq m \leq (k+1)j + k$, $j \in \mc J_n$
  we have
  \begin{equation*}
    \Efm\{ W_m^2 \} = (\Var\{ W \})_{m, m} = (\Var\{ \gamma_{I(x_j,
      h)} \})_{m - (k+1)j, m - (k+1)j} = 1,
  \end{equation*}
  since $\Var\{ \gamma_{I(x_j, h)} \} = \mb I_{k+1}$. Then, we have
  proved that on $\cap_{j \in \mc J_n} \Omega_{n, I(x_j, h)}$,
  \begin{align*}
    \max_{j \in \mc J_n} \norminfty{ \mb \Lba_{I(x_j, h)}^{-1} ( \wh
      \tta_{I(x_j, h)} &- \tta_{I(x_j, h)})} \\
    &\leq \lba^{-1} ( \mc G_{I(x_j, h)}) (k+1) \big( L h^s +
    \frac{\sigma}{\sqrt{n \bar \mu_n(I(x_j, h))}} W^M \big),
  \end{align*}
  where $W^M = \max_{0 \leq m \leq (k+1) |\mc J_n|} |W_{m}|$. Since
  $\mc C_n \subset \mrm N_n \cap \Omega_n \cap \mc L_n$, we have on
  $\mc C_n$ for $h = h_n$ or $h = t_n$,
  \begin{align*}
    \max_{j \in \mc J_n} \norminfty{ \mb \Lba_{I(x_j, h)}^{-1} ( \wh
      \tta_{I(x_j, h)} &- \tta_{I(x_j, h)})} \\
    &\leq (1 + o(1)) \lba^{-1}(\mc G) (k+1) \big( L h^s +
    \frac{\sigma}{\sqrt{n h \mu_j}} W^M \big).
  \end{align*}
  Since $\mc C_n \subset \mc D_n$, we have for any $j \in \mc J_n$, $0
  \leq m \leq k$,
  \begin{equation*}
    \mc C_n \subset \bar{\mrm D}_{n, 2m, I(x_j, h_n), \delta_n} \cap
    \bar {\mrm D}_{n,2m, I(x_j, t_n), \delta_n},
  \end{equation*}
  thus on $\mc C_n$, when $h = h_n$ or $h = t_n$, we clearly have
  \begin{equation*}
    (\mb \Lba_{I(x_j, h)})_{m,m} = \norm{\phi_{I(x_j, h), m}}_{I(x_j,
      h)}^{-1} \leq (1 + o(1)) h^{-m} \sqrt{2m + 1}.
  \end{equation*}
  Since $\wt f_n^{(m)}(x_j) - f^{(m)}(x_j) = m! \big( ( \wh
  \tta_{I(x_j, h_n)})_m - (\tta_{I(x_j, h_n)})_m \big)$, it follows
  that on $\mc C_n$:
  \begin{align*}
    |\wt f_n^{(m)}(x_j) &- f^{(m)}(x_j)| \\
    &\leq (1 + o(1)) \lba^{-1}(\mc G) m! \sqrt{2m + 1} (k+1) h_n^{-m}
    (L h_n^s + \frac{\sigma}{\sqrt{n h_n \mu_j}} W^M) \\
    &\leq (1 + o(1)) C L h_n^{s-m} (1 + (\log n)^{-1/2} W^M),
  \end{align*}
  thus~\eqref{eq:bias_variance_derivatives_estimates}.
  Inequality~\eqref{eq:bias_variance_boundary_estimator} is obtained
  similarly.
\end{proof}

\begin{proof}[Proof of lemma~\ref{lem:LPA_bad_case}]
  \label{proof:lemma_bad_case_chap3}
  If $\bar \mu_n(I) = 0$ we have $\wh \tta_I = 0$ and the result is
  obvious, thus we assume $\bar \mu_n(I) > 0$. In this case, $\mb
  \Lba_I$, $\mb{\bar X}_I$ and $\mc G_I$ are invertible, and by
  definition of $\wh \tta_I$,
  \begin{equation*}
    \wh \tta_I = \mb \Lba_I \mb \Lba_I^{-1} \wh \tta_I = \mb \Lba_I
    \mc G_I^{-1} \mb \Lba_I \bar{\mb X}_I \wh \tta_I = \mb \Lba_I \mc
    G_I^{-1} \mb \Lba_I \mb Y_I = \mb \Lba_I \mc G_I^{-1} (\mb B_I + \mb
    V_I),
  \end{equation*}
  where $(\mb B_I)_m = \norm{\phi_{I, m}}_I^{-1} \prodsca{f}{\phi_{I,
      m}}_I$ and $(\mb V_I)_m = \norm{\phi_{I, m}}_I^{-1}
  \prodsca{\xi}{\phi_{I, m}}_I$. Since $\norminfty{f} \leq Q$ we have
  $|(\mb B_I)_m| = \norm{\phi_{I, m}}_I^{-1} | \prodsca{f}{\phi_{I,
      m}}_I | \leq \norm{f}_I \leq Q$, thus $\norm{\mb B_I}_{\infty}
  \leq Q$.

  Conditionally on $\mf X_n$, $\mb V_I$ is centered Gaussian and it is
  an easy computation to see that its covariance matrix is equal to
  $\sigma^2 (n \bar \mu_n(I))^{-1} \mb \Lba_I \mb X_I \mb \Lba_I$.
  Then $\mb \Lba_I \mc G_I^{-1} \mb V_I$ is conditionally on $\mf X_n$
  centered Gaussian with covariance matrix $\sigma^{2} (n \bar
  \mu_n(I))^{-1} \bar{\mb X}_I^{-1} \mb X_I \bar{\mb X}_I^{-1}$. If
  $e_m$ is the canonical vector with coordinates $(e_m)_p = \ind{p =
    m}$, we have
  \begin{equation*}
    | (\wh \tta_{I})_m | = |\prodsca{\wh \tta_I}{e_m}| = |\prodsca{\mb
      \Lba_I \mc G_I^{-1} \mb B_I}{e_m} | + \sigma \sqrt{k+1} \,
    \gamma,
  \end{equation*}
  where $\gamma = (\sigma \sqrt{k+1})^{-1} \prodsca{\mb \Lba_I \mc
    G_I^{-1} \mb V_I}{e_m}$. By definition, we have $\norm{\bar{\mb
      X}_I^{-1}} = \lba^{-1}(\bar{\mb X}_I) \leq \sqrt{n \bar
    \mu_n(I)}$, and clearly $\norm{\mb X_I} \leq k + 1$ and $\norm{\mb
    \Lba_I^{-1}} \leq 1$. Then, conditional on $\mf X_n$, $\gamma$ is
  centered Gaussian with variance
  \begin{equation*}
    \frac{ \prodsca{e_m}{\bar{\mb X}_I^{-1} \mb X_I \bar{\mb X}_I^{-1}
        e_m} }{(k+1) n \bar \mu_n(I)}  \leq \frac{\norm{\bar{\mb
          X}_I^{-1}}^2 \norm{\mb X_I}}{(k+1) n \bar \mu_n(I)} \leq 1.
  \end{equation*}
  Since $\norm{\mc G_I^{-1}} \leq \norm{\mb \Lba_I^{-1}}
  \norm{\bar{\mb X}_{I}^{-1}} \norm{\mb \Lba_I^{-1}} \leq \sqrt{n \bar
    \mu_n(I)} \leq \sqrt{n}$ and $(\mb \Lba_I)_{0, 0} = 1$, we have
  \begin{equation*}
    \Efm\big\{ |(\wh \tta_I)_0|^p | \mf X_n \big\} \leq (k+1)^{p/2}
    n^{p / 2} (Q \vee 1)^p \Efm\{\big (1 + \sigma |\gamma|)^p | \mf
    X_n \big\} = O(n^{p / 2}),
  \end{equation*}
  for any $I \subset [0, 1]$, and since $\norm{\mb \Lba_I} \leq
  \sqrt{n}$ on $\Gamma_{n, I}$, it follows that
  \begin{equation*}
    \Efm\big\{ |(\wh \tta_I)_m|^p | \mf X_n \big\} \leq (k+1)^{p/2}
    n^p (Q \vee 1)^p \Efm\{\big (1 + \sigma |\gamma|)^p | \mf X_n
    \big\} = O(n^p),
  \end{equation*}
  for any $1 \leq m \leq k$.
\end{proof}

\begin{proof}[Proof of lemma~\ref{lem:risk_expectations_bounded}]
  We show that for any $p > 0$,
  \begin{equation}
    \label{eq:prop_expect_bounded_proof_1}
    \sup_{f \in \Sigma^Q(s, L)} \Efm\{ \mc E_{n, f}^p \} =
    O(n^{ p(1 + s/(2s + 1)) }),
  \end{equation}
  which entails~\eqref{eq:expectation_poynomially_bounded}. By
  definition of $H_n(x)$, we have $H_n(x) \geq (\log n / n)^{1 /
    (2s)}$ for any $x \in [0, 1]$.  Since $\norminfty{f} \leq Q$, we
  have for any $j \in \mc J_{2, n}$,
  \begin{equation*}
    |\wh f_n(x_j)| \leq \delta_n^{-1} (n / \log n)^{1 / (2s)} (Q + |
    \bar \xi_n| / \sqrt{n}) \norminfty{K_s},
  \end{equation*}
  where $\bar \xi_n = \sumin \xi_i / \sqrt{n}$ is standard Gaussian.
  Then,
  \begin{align*}
    \Efm\big\{ |\wh f_n(x_j)|^p | \mf X_n \big\} &\leq \delta_n^{-p} (
    (n / \log n)^{p / (2s)} (Q \vee 1)^p \Efm\{ (1 + |\bar \xi_n| )^p
    | \mf X_n \} \norminfty{K_s} \\ &= O(n^{p / (2 s)} (\log n)^{p(1 -
      1/ (2s))}).
  \end{align*}
  When $j \in \mc J_{n, 1} \cup \mc J_{n, 3}$, we have $\wh f_n(x_j) =
  \wh \tta_{I(x_j, t_n)}$ and in view of lemma~\ref{lem:LPA_bad_case},
  \begin{equation*}
    \Efm\big\{ |\wh f_n(x_j)|^p | \mf X_n \big\}= O(n^{p / 2}).
  \end{equation*}
  For any $j \in \mc J_n$, since $\wt f_n^{(m)}(x_j) = m! ( \wh
  \tta_{I(x_j, h_n)} )_m$, we have in view of
  lemma~\ref{lem:LPA_bad_case} that on $\Gamma_{n, I(x_j, h_n)}$,
  \begin{equation*}
    \Efm\big\{ |\wt f_n^{(m)}(x_j)|^p | \mf X_n \big\} = O(n^p),
  \end{equation*}
  for any $1 \leq m \leq k$. Then, we obtain that for any
  $\norminfty{f} \leq Q$,
  \begin{align*}
    \mc E_{n, f} = O((n / \log n)^{s / (2s + 1)}) \big( \sup_{x \in
      [0, 1]} |\wh f_n(x)| + Q \big),
  \end{align*}
  and since
  \begin{equation*}
    \sup_{x \in [0, 1]} |\wh f_n(x)| \leq \max_{j \in \mc J_n} \Big( |
    \wh f_n(x_j)| + \big( \sum_{m = 1}^k \frac{| \wt f_n^{(m)}(x_j) |
    }{m!} \big) \ind{\Gamma_{n, I(x_j, h_n)}} \Big) = O(n^p),
  \end{equation*}
  thus~\eqref{eq:prop_expect_bounded_proof_1}
  and~\eqref{eq:expectation_poynomially_bounded}.
\end{proof}

\begin{proof}[Proof of lemma~\ref{lem:awful_lemma}]
  \label{proof:awful_lemma_chap3}
  The proof is divided in several steps. We recall that $q_j = n c_s
  h_j \mu_j$ and $\bar q_j = n c_s H_j \mu_j$. \\
  
  \noindent \emph{Step 1.} We prove that for any $j \in \mc J_{2, n}$
  and $n$ large enough,
  \begin{equation}
    \label{eq:lemma_deviations_1}
    \Pm\{ \mrm B_{n, j}^c \} \leq 2 \exp( -D_1 \delta_n^2 n^{2s / (2s
      + 1)}),
  \end{equation}
  where $D_1$ is a positive constant. Consider the sequence of i.i.d
  variables $\zeta_{i, j} \eqdef K_{i,j} - \Em\{ K_{i,j} \}$, $1 \leq
  i \leq n$.
  Since $\mu \in \Sigma_q(\nu, \varrho)$ and $\int K = 1$, we have for
  $n$ large enough $| \Em\{ K_{1, j}\} / q_j - 1 | \leq \delta_n / 2$,
  thus $\mrm B_{n, j}^c \subset \big\{ | \sumin \zeta_{i, j} | / q_j
  \leq \delta_n / 2 \big\}.$ Since $|\zeta_{i,j}| \leq 2
  \norm{K}_{\infty}$ and for $n$ large enough $\sumin \Em\{
  \zeta_{i,j}^2 \} \leq (1 + \delta_n) q_j \int K^2$, the Bernstein
  inequality entails \eqref{eq:lemma_deviations_1}. \\
  
  \noindent \emph{Step 2.} We prove that for any $j \in \mc J_{n, 2}$,
  \begin{equation}
    \label{eq:lemma_deviations_2}
    \Pm\{ \mrm A_{n, j}^c \cap \mrm C_{n,j} \} \leq 2 \exp( -D_2
    \delta_{2, n}^2 \, n^{2s / (2s + 1)}),
  \end{equation}
  where $D_2$ is a positive constant and $\delta_{2,n} \eqdef
  \delta_n^{s_1}$, $s_1 = s \wedge 1$. In view of
  \eqref{eq:use_of_K_regularity}, we have on $\mrm C_{n,j}$
  \begin{equation}
    \label{eq:awful_lemma_proof_1}
    | \bar K_{i,j} - K_{i,j} | \leq \kpa T_s^{s_1} \Big(
    \frac{\delta_n}{1 - \delta_n} \Big)^{s_1} \ind{\mrm M_{i,j}}
  \end{equation}
  where we recall that $\mrm M_{i,j} = \{ |X_i - x_j| \leq c_s T_s (1
  + \delta_n) h_j\}$. We define $\eta_{i,j} \eqdef \ind{\mrm M_{i,j}}
  - \Pm\{ \mrm M_{i,j} \}$. On $\mrm C_{n,j}$ we have for $n$ large
  enough $2 c_s T_s H_n^M \leq \delta_n$, and since $x_j \in [\tau_n,
  1 - \tau_n]$,
  \begin{align*}
    x_j \leq 1 - \tau_n = 1 - 2 c_s T_s H_n^M &\leq 1 - 2 c_s T_s H_j
    \\
    &\leq 1 - 2 c_s T_s (1 - \delta_n) h_j \leq 1 - c_s T_s (1 +
    \delta_n) h_j
  \end{align*}
  for $n$ large enough.  On the other hand we have similarly $x_j \geq
  c_s T_s (1 + \delta_n) h_j$. Thus, since $\mu \in \Sigma_q(\nu,
  \varrho)$ we have
  \begin{equation}
    \label{eq:awful_lemma_proof_2}
    \Big| \frac{\Pm\{ \mrm M_{i,j} \}}{(1 + \delta_n) c_s h_j \mu_j} -
    2 T_s \Big| \leq \frac{1}{q} \int_{|y| \leq T} |\mu(x_j + c_s y (1
    + \delta_n) h_j) - \mu_j| dy = O(h_n^{\nu}).
  \end{equation}
  Since $x_j \in [c_s T_s (1 + \delta_n) h_j, 1 - (1 + \delta_n) c_s
  T_s h_j] \subset [ c_s T_s h_j, 1 - c_s T_s h_j]$, we have for $n$
  large enough on $\mrm C_{n,j}$,
  \begin{equation}
    \label{eq:awful_lemma_proof_3}
    \begin{split}
      \Big| \frac{ \Efm\{ K_{1,j} \} }{c_s H_j \mu_j} - 1 \Big| &\leq
      \frac{h_j}{H_j \mu_j} \int |K(y)| |\mu(x_j + y c_s h_j) -
      \mu_j| dy + \Big| \frac{h_j}{H_j} - 1 \Big| \\
      &\leq O(h_n^{\nu}) + \frac{\delta_n}{1 - \delta_n}.
    \end{split}
  \end{equation}
  Then, combining~\eqref{eq:awful_lemma_proof_1},
  \eqref{eq:awful_lemma_proof_2} and~\eqref{eq:awful_lemma_proof_3} we
  obtain that on $\mrm C_{n, j}$ and for $n$ large enough,
  \begin{align*}
    \Big| \frac{\sumin \bar K_{i,j}}{\bar q_j} - 1 \Big| &\leq
    \frac{o(1)}{\bar q_j} | \sumin \eta_{i,j} | + \frac{ \kpa
      T^{s_1}\delta_n^{s_1}}{ (1 - \delta_n)^{s_1}} \frac{\Pm\{ \mrm
      M_{1, j} \}}{c_s H_j \mu_j} \\
    &+ \frac{1}{\bar q_j} | \sumin \zeta_{i,j} | + O(h_n^{\nu}) +
    \frac{\delta_n}{1 - \delta_n} \\
    &\leq \frac{o(1)}{q_j} | \sumin \eta_{i,j} | + \frac{1 +
      o(1)}{q_j} | \sumin \zeta_{i,j} | + 2(2 \kpa T^{s_1 + 1} + 1)
    \delta_n^{s_1},
  \end{align*}
  and taking $L_1 \eqdef 4(\kpa T^{s_1 + 1} + 1)$, we obtain
  \begin{equation*}
    \Pfm\{ \mrm A_{n,j}^c \cap \mrm C_{n,j} \} \leq \Pm\big\{| \sumin
    \eta_{i,j} | > \delta_n^{s_1} q_j \big\} + \Pm\{ | \sumin
    \zeta_{i,j} | > \delta_n^{s_1} q_j / 2\}.
  \end{equation*}
  Then, applying Bernstein inequality to the sum of variables
  $\eta_{i,j}$ and $\zeta_{i,j}$, $1 \leq i \leq n$, we
  obtain~\eqref{eq:lemma_deviations_2}. We can prove
  \begin{equation}
    \label{eq:lemma_deviation_4}
    \Pm\{ \mrm E_{n, j}^c \cap \mrm C_{n,j} \} \leq 2 \exp( -D_3
    \delta_{2, n}^2 n^{2s / (2s + 1)}),
  \end{equation}
  where $D_3$ is a positive constant in the same way as for the proof
  of~\eqref{eq:lemma_deviations_2} with a good choice for $L_2$. \\

  \noindent \emph{Step 3.} We define the event
  \begin{equation*}
    \mrm D_{n, m, I(x, h), \delta} \eqdef \bigg\{ \Big|
    \frac{1}{\mu(x) h^{m+1}} \int_{I(x, h)} \phi_{I(x, h), m} \, d\bar
    \mu_n - \chi_m \Big| \leq \delta \bigg\},
  \end{equation*}
  and we prove that if $\delta_{1,n} \eqdef 1 - (1 + \delta_n)^{-(2s +
    1)}$,
  \begin{equation}
    \label{eq:lemma_deviation_3}
    \mrm D_{n, 0, I(x_j, (1 - \delta_n) h_j), \delta_{1, n}} \cap \mrm
    D_{n, 0, I(x_j, (1 + \delta_n) h_j), \delta_{1, n}} \subset \mrm
    C_{n,j}.
  \end{equation}
  From the definitions of $H_j$ and $h_j$ (see section
  \ref{sec:confidence_bands}) we obtain
  \begin{align*}
    \{ (1 - \delta_n) h_j < H_j \} &= \big\{ (1 - \delta_n)^{2s}
    h_j^{2s} < \log n / \big( n \bar \mu_n(I(x_j, (1 - \delta_n)
    h_j)) \big) \big\} \\
    &= \Big\{ \frac{\bar \mu_n( I(x_j, (1 - \delta_n) h_j) )}{\mu_j (1
      - \delta_n) h_j} \leq (1 - \delta_n)^{-(2s+1)} \Big\},
  \end{align*}
  and then
  \begin{equation*}
    \mrm D_{n, 0, I(x_j, (1 - \delta_n) h_j), \delta_{1, n}} \subset
    \{ (1 - \delta_n) h_j < H_j \}.
  \end{equation*}
  We can prove in the same way that on the other hand,
  \begin{equation*}
    \mrm D_{n, 0, I(x_j, (1 + \delta_n) h_j), \delta_{1, n}} \subset
    \{ (1 + \delta_n) h_j \geq H_j \},
  \end{equation*}
  hence~\eqref{eq:lemma_deviation_3}. \\
  
  \noindent \emph{Step 4.} We
  prove~\eqref{eq:mc_A_n_subset_many_events}. If $\delta_{3, n} =
  \delta_n / (2 - \delta_n)$, we clearly have for any interval~$I$,
  \begin{equation*}
    \mrm D_{n, m, I, \delta_{3, n}} \cap \mrm D_{n, 0, I,
      \delta_{3, n}} \subset \bar{\mrm{D}}_{n, m, I, \delta_{n}}.
  \end{equation*}
  Using the fact that $\lba(M) = \inf_{\norm{x} = 1} \prodsca{x}{M x}$
  for any symmetrical matrix $M$ and since $\mc G_{I}$, $\mc G$, $\mb
  X_I$ are symmetrical, it is easy to see that
  \begin{equation}
    \label{eq:awful_lemma_proof_subset_mc_L_n}
    \bigcap_{0 \leq p, q \leq k} \Big\{ | (\mc G_{I} - \mc G)_{p,q} |
    \leq \frac{\delta_{n}}{(k+1)^2} \Big\} \subset \mc L_{n, I},
  \end{equation}
  and that
  \begin{align*}
    \bigcap_{m = 0} ^{2k} \bar{\mrm{D}}_{n, m, I, \frac{\delta_n}{
        (k+1)^2 }} &\subset \bigcap_{0 \leq p, q \leq k} \Big\{ \big|
    (\mb X_{I} - \mb X)_{p,q} \big| \leq \frac{\delta_n}{(k+1)^2}
    \Big\} \\
    &\subset \mc \{ | \lba(\mb X_{I}) - \lba(\mb X) | \leq \delta_{n}
    \}.
  \end{align*}
  Recalling that if $I = I(x_j, h)$,
  \begin{equation*}
    (\mc G_{I})_{p, q} = \frac{\prodsca{\phi_{I, p}}{\phi_{I, q}}_I}{
      \norm{ \phi_{I, p}}_I \norm{\phi_{I, q}}_I} = \frac{ \frac{1}{
        \mu_j h^{m+1}} \int_{I} \phi_{I, p+q}\, d \bar
      \mu_n }{ \sqrt{\frac{1}{\mu_j h^{m+1}} \int_{I} \phi_{I,
          2p}\, d\bar \mu_n}  \sqrt{\frac{1}{\mu_j h^{m+1}}
        \int_{I} \phi_{I, 2q}\, d\bar \mu_n} } ,
  \end{equation*}
  it is easy to see that if $\delta_{4, n} = \delta_n / \big( (2 -
  \delta_n)(2k + 1)(k+1)^2 \big)$,
  \begin{equation*}
    \mrm D_{n, 2 p, I, \delta_{4, n}} \cap \mrm D_{n, 2 q, I,
      \delta_{4, n}} \cap \mrm D_{n, p + q, I, \delta_{4, n}}
    \subset \Big\{ | (\mc G_{I} - \mc G)_{p,q} |\leq \frac{ \delta_{n}
    }{(k+1)^2} \Big\},
  \end{equation*}
  thus
  \begin{equation*}
    \bigcap_{m = 0}^{2k} \mrm D_{n, m, I, \delta_{4, n}} \subset
    \mc L_{n, I},
  \end{equation*}
  and clearly for $n$ large enough, if $I = I(x_j, h_n)$ or $I =
  I(x_j, t_n)$,
  \begin{equation}
    \label{eq:awful_lemma_proof_subset_Omega}
    \bigcap_{m = 0}^{2k} \mrm D_{n, m, I, \delta_{4, n}} \subset
    \{ | \lba(\mb X_I) - \lba(\mb X) | \leq \delta_n \} \cap \Big\{
    \Big| \frac{ \bar \mu_n(I) }{|I| \mu_j} - 1 \Big| \leq \delta_n
    \Big\} \subset \Omega_{n, I}.
  \end{equation}
  Moreover, if $I = I(x_j, h_n)$, we have on $\bar {\mrm D}_{n, 2m, I,
    \delta_n}$ for any $1 \leq m \leq k$ and $n$ large enough,
  \begin{equation}
    \label{eq:bar_D_n_2m_subset_Gamma_I}
    \norm{\phi_{I, m}}_I \geq (1 - o(1)) h_n^m \sqrt{2m + 1} \geq 1 /
    \sqrt{n}.
  \end{equation}
  We define
  \begin{multline*}
    \mrm D_{n, m} \eqdef \bigcap_{j \in \mc J_n} \Big( \mrm D_{n, m,
      I(x_j, h_n), \delta_{5, n}} \cap \mrm D_{n, m, I(x_j, t_n),
      \delta_{5, n}} \\ \cap \mrm D_{n, 0, I(x_j, (1 - \delta_n) h_j),
      \delta_{5, n}} \cap \mrm D_{n, 0, I(x_j, (1 + \delta_n) h_j),
      \delta_{5, n}} \Big),
  \end{multline*}
  where $\delta_{5, n} = \delta_{4, n} \wedge \delta_{3, n} \wedge
  \delta_{1, n}$, $\mrm D_n = \bigcap_{m=0}^{2k} \mrm D_{n, m}$ and we
  choose
  \begin{equation*}
    \mc A_n \eqdef \mrm D_{n} \cap \mrm A_{n} \cap \mrm B_{n} \cap
    \mrm E_{n}.
  \end{equation*}
  In view of~\eqref{eq:lemma_deviation_3},
  \eqref{eq:awful_lemma_proof_subset_mc_L_n},
  \eqref{eq:awful_lemma_proof_subset_Omega},
  \eqref{eq:bar_D_n_2m_subset_Gamma_I} we have $\mc A_n \subset \mrm
  C_n \cap \Omega_n \cap \mc L_n \cap \Gamma_n$ and since $\mrm D_{n,
    0, I, \delta} = \mrm N_{n, I}$ we
  obtain~\eqref{eq:mc_A_n_subset_many_events}. \\

  \noindent \emph{Step 5.} We prove~\eqref{eq:mc_A_n_deviation}.
  Using Bernstein inequality, it is easy to show that for $n$ large
  enough, if $h = h_n$, $h = t_n$, $h = (1 - \delta_n) h_j$ or $h = (1
  + \delta_n) h_j$,
  \begin{equation*}
    \Pm\{ \mrm D_{n, m, I(x_j, h), \delta_{5, n}}^c \} \leq 2 \exp(
    -D_4 \delta_{5, n}^2 n h) \leq 2 \exp( -D_5 n^{ s / (2s + 1) }),
  \end{equation*}
  with $D_4, D_5$ positive constants, where we used the fact that
  $\delta_{5, n}^2 n^{s / (2s + 1)} > 1$ for $n$ large enough and $n h
  \geq D_6 n^{2s / (2s + 1)}$. In view of~\eqref{eq:lemma_deviation_3}
  we have $\mrm D_n \subset \mrm C_n$, hence
  \begin{align*}
    \Pm\{ \mc A_n^c \} &\leq \Pfm\{ \mrm D_n^c \} + \Pfm\{ \mrm A_n^c
    \cap \mrm C_n \} + \Pfm\{ \mrm B_n^c \cap \mrm C_n \} \\
    &+ \Pfm\{ \mrm E_n^c \cap \mrm C_n \} + 3 \Pfm\{ \mrm C_n^c \} \\
    &\leq 4\, \Pfm\{ \mrm D_n^c \} + \Pfm\{ \mrm A_n^c \cap \mrm C_n
    \} + \Pfm\{ \mrm B_n^c \cap \mrm C_n \} + \Pfm\{ \mrm E_n^c \cap
    \mrm C_n \} \\
    &\leq 2 (8k + 7) M_n \exp( - 2 D_{\mc A} n^{s / (2s + 1)}) \leq
    \exp( -D_{\mc A} n^{s / (2s + 1)} ),
  \end{align*}
  for $n$ large enough, where $D_{\mc A} \eqdef (D_1 \vee D_2 \vee D_3
  \vee D_5) / 2$, where we used~\eqref{eq:lemma_deviations_1},
  \eqref{eq:lemma_deviations_2} and \eqref{eq:lemma_deviation_4}.
\end{proof}

\section{Proof of theorem~2}
\label{sec:proof_of_the_lower_bound}

The proof of the lower bound is heavily based on arguments found in
\cite{korostelev93}, \cite{donoho94}, \cite{korostelev_nussbaum_99}
and \cite{bertin02}. It is mainly a modification of the former proof
in \cite{bertin02}. It consists in a classical reduction to the
Bayesian risk over an hardest cubical subfamily of functions, see for
instance \cite{donoho94}. The main difference with the former proofs
is that the subfamily of functions depends on the design via the
bandwidth $h_{n, \mu}(x)$, which is adapted to the local amount of
data.

\subsection{Preparatory results}
\label{sec:lower_bound_preparatory_results}

We begin with some definitions. We recall that $\varphi_s$ is defined
by~\eqref{eq:OR_optimal_function} and that it has a compact support
$[-T_s, T_s]$. Let $h_n^I \eqdef \max_{x \in I_n} h_{n, \mu}(x)$ and
\begin{equation*}
  \Xi_n = 2 T_s c_s (2^{1 / (s - k)} + 1) h_n^I.
\end{equation*}
If $I_n = [a_n, b_n]$, $M_n = [|I_n| \, \Xi_n^{-1}]$, we define the
points
\begin{equation}
  \label{eq:lower_bound_points}
  x_j = a_n + j \, \Xi_n, \quad j \in \mc J_n \eqdef \{ 1, \ldots, M_n
  \}.
\end{equation}
In order to unload the notations, we denote again $\mu_j = \mu(x_j)$,
$h_j = h_{n, \mu}(x_j)$.

\begin{lemma}
  \label{lem:H_n_deviation}
  Let define the event
  \begin{equation*}
    \mrm H_{n, j} \eqdef \Big\{ \Big| \frac{1}{n c_s h_j \mu_j} \sumin
    \varphi_s^2 \Big( \frac{X_i - x_j}{c_s h_j} \Big)- 1 \Big|
    \leq \von \Big\},
  \end{equation*}
  and $\mrm H_n \eqdef \cap_{j \in \mc J_n} \mrm H_{n, j}$. We have
  \begin{equation*}
    \lim_{n \raro +\infty} \Pm\{ \mrm H_n \} = 1.
  \end{equation*}
\end{lemma}

\begin{proof}
  We use Bernstein inequality to the sum of variables
  $\varphi_s^2((X_i - x_j) / (c_s h_j))$, for $1 \leq i \leq n$, where
  we use the fact that $\norm{\varphi_s}_2 = 1$ (see
  section~\ref{sec:optimal_recovery}) and we derive a deviation
  inequality for the events $\mrm H_{n, j}^c$. Then, bounding from
  above the probability of $\cup_{j \in \mc J_n} \mrm H_{n, j}^c$ by
  the probabilities sum, the result follows easily.
\end{proof}

The subfamily of functions is defined as follows. We consider an
hypercube $\Theta \subset [-1, 1]^{M_n}$, and for $\tta \in \Theta$ we
define the functions
\begin{equation*}
  f(x; \tta) = \sum_{j \in \mc J_n} \tta_j f_j(x), \quad f_j(x) = L
  c_s^s h_j^s \varphi_s \Big( \frac{x - x_j}{c_s h_j} \Big).
\end{equation*}
Clearly, $f_j \in \Sigma(s, L)$. Let us show that $f(\cdot \,; \tta)
\in \Sigma(s, L)$. We note that
\begin{equation*}
  \supp\Big( \varphi_s \big( \frac{\cdot - x_j}{c_s h_j} \big) \Big)
  = \big[x_j - c_s T_s h_j, \, x_j + c_s T_s h_j \big] \eqdef I_j.
\end{equation*}
If $x, y \in I_j$ then $f(x; \tta) = \tta_j f_j(x)$, $f(y; \tta) =
\tta_{j} f_{j}(y)$ and the result is obvious. To complete the proof,
it suffices to consider the case $x \in I_j$ and $y \in I_{j+1}$. In
this case, we have
\begin{align*}
  |f^{(k)}(x;& \tta) - f^{(k)}(y; \tta)| \\
  &= | \tta_j f_j^{(k)}(x) - \tta_{j+1} f_{j+1}^{(k)}(y)| \\
  &\leq |f_{j}^{(k)}(x) - f_j^{(k)}(x_j + c_s T_s h_j)| +
  |f_{j+1}^{(k)}(x_{j+1} - c_s T_s h_{j+1}) - f_{j+1}^{(k)}(y)| \\
  &\leq L \big( |x - x_j - c_s T_s h_j|^{s-k} + |x_{j+1} - c_s T_s
  h_{j+1} -y |^{s - k} \big) \\
  &\leq L \big( (2 c_s T_s h_j)^{s - k} + (2 c_s T_s h_{j+1})^{s - k}
  \big) \leq 2 L (2 c_s T_s h_n^I )^{s-k}.
\end{align*}
Moreover, since $x \in I_j$ and $y \in I_{j+1}$ we have
\begin{equation*}
  |x - y| \geq x_{j+1} - x_j - c_s T_s (h_j +h_{j+1}) \geq \Xi_n - 2
  c_s T_s h_n^I = 2^{1 / (s - k)} (2 c_s T_s h_n^I),
\end{equation*}
and finally
\begin{equation}
  \label{eq:proof_f_x_tta_belongs_to_sigma}
  |f^{(k)}(x; \tta) - f^{(k)}(y; \tta)| \leq L |x - y|^{s - k},
\end{equation}
thus $f(\cdot\,; \tta) \in \Sigma(s, L)$. For any $j \in \mc J_n$, we
define the statistics
\begin{equation*}
  y_j = \frac{\sumin Y_i \varphi_s(X_i)}{\sumin \varphi_s^2(X_i)}.
\end{equation*}

\begin{lemma}
  \label{lem:lemma_y_j_sufficiency}
  Conditionally on $\mf X_n$, the $y_j$ are Gaussian and independent.
  Moreover, if $v_j^2 = \Efm\{ y_j^2 | \mf X_n \}$, we have on $\mrm
  H_{n, j}$
  \begin{equation}
    \label{eq:control_variance_y_j}
    \Efm\{ y_j | \mf X_n \} = \tta_j, \quad \frac{2s + 1}{2 (1 + \von)
      \log n} \leq v_j^2 \leq \frac{2s + 1}{2 (1 - \von) \log n}.
  \end{equation}
  In the model~\eqref{eq:model_regression} with $f(\cdot) = f(\cdot\,;
  \tta)$, conditionally on $\mf X_n$, the likelihood function of
  $(Y_1, \ldots , Y_n)$ can be written on $\mrm H_n$ in the form
  \begin{equation*}
    \frac{\mrm d\Pfm}{\mrm d \lba^n}|_{\mf X_n}(Y_1, \ldots, Y_n) =
    \prod_{i=1}^n g_{\sigma}(Y_i) \prod_{j \in \mc J_n}
    \frac{g_{v_j}(y_j - \tta_j)}{g_{v_j}(y_j)},
  \end{equation*}
  where $g_{v}$ is the density of $\mc N(0, v^2)$, and $\lba^n$ is the
  Lebesgue measure over $\setR^n$.
\end{lemma}

\begin{proof}
  By construction the $f_j$ have disjoint supports, thus it is easy to
  see that conditionally on $\mf X_n$ the $y_j$ are Gaussian
  independent with conditional mean $\tta_j$. Using the definition of
  $\mrm H_n$ and since
  \begin{equation*}
    \Efm\{ y_j^2 | \mf X_n \} = \frac{\sigma^2}{\sumin f_j^2(X_i)},
  \end{equation*}
  it is an easy computation to see that on $\mrm H_n$, we
  have~\eqref{eq:control_variance_y_j}. The last part of the
  lemma follows from the following computation: 
  \begin{align*}
    \prod_{i=1}^n &g_{\sigma}(Y_i) \prod_{j \in \mc J_n}
    \frac{g_{v_j}(y_j - \tta_j) }{g_{v_j}(y_j)} \\
    &= \frac{1}{\sigma^n (2 \pi)^{n/2}} \prod_{i=1}^n \exp\big( -Y_i^2
    / (2 \sigma^2) \big) \prod_{j \in \mc J_n} \exp\big( (2 \tta_j y_j
    - \tta_j) / (2 v_j^2) \big) \\
    &= \frac{1}{\sigma^n (2 \pi)^{n/2}} \prod_{i=1}^n \bigg[ \exp\Big(
    \frac{ - Y_i^2 + \sum_{j \in \mc J_n} \big( 2 Y_j \tta_j f_j(X_i)
      - \tta_j^2 f_j(X_i)^2 \big)}{2 \sigma^2} \Big) \bigg] \\
    &= \frac{1}{\sigma^n (2 \pi)^{n/2}} \prod_{i=1}^n \exp\Big( -
    \frac{(Y_i - f(X_i; \tta))^2}{2 \sigma^2} \Big) = \frac{\mrm d
      \Pfm}{\mrm d\lba^n} |_{\mf X_n} (Y_1, \ldots, Y_n). \qedhere
  \end{align*}
\end{proof}

\subsection{Proof of theorem \ref{thm:lower_bound}}
\label{sec:proof_of_lower_bound}

We denote in the following $\Sigma = \Sigma(s, L)$ and $\mc E_{n, f,
  T}^I = \sup_{x \in I} r_{n, \mu}(x)^{-1} |T(x) - f(x)|$. Since $w$
is nondecreasing and $f(\cdot\, ; \tta) \in \Sigma$ for any $\tta \in
\Theta$, we have for any distribution $\mc B$ on $\Theta$ by a
minoration of the minimax risk by the Bayesian risk,
\begin{align*}
  \inf_{T} \sup_{f \in \Sigma} \Efm\big\{ w( \mc E_{n, f, T}^I)\big\}
  &\geq w\big((1 - \von) P\big) \inf_{T} \sup_{f \in
    \Sigma} \Pfm\big\{ \mc E_{n, f, T}^I \geq (1 - \von) P \big\} \\
  &\geq w\big((1 - \von) P\big) \inf_{T} \int_{\Theta} \Prob_{\tta}^n
  \big\{ \mc E_{n, f, T}^I \geq (1 - \von) P \big\} \mc B(d \tta),
\end{align*}
where $\Prob_{\tta}^n = \Prob_{f(\cdot\,; \tta), \mu}^n$. Since by
construction $f(x_j; \tta) = r_j \tta_j P$ and $x_j \in I_n$, we
obtain
\begin{align*}
  \inf_{T} \int_{\Theta} \Prob_{\tta}^n &\big\{ \mc E_{n, f, T}^I \geq
  (1 - \von) P \big\} \mc B(d \tta) \\
  &\geq \inf_{\wh \tta} \int_{\Theta} \int_{\mrm H_n} \Prob_{\tta}^n
  \big\{ \max_{j \in \mc J_n} |\wh \tta_j - \tta_j
  | \geq 1 - \von | \mf X_n \big\} \mrm d\Pm \mc B(d\tta), \\
  &\geq \int_{\mrm H_n} \inf_{\wh \tta} \int_{\Theta} \Prob_{\tta}^n
  \big\{ \max_{j \in \mc J_n} |\wh \tta_j - \tta_j | \geq 1 - \von |
  \mf X_n \big\} \mc B(d\tta) \mrm d\Pm,
\end{align*}
where $\inf_{\wh \tta}$ is taken among any measurable vector (with
respect to the observations~\eqref{eq:model_regression}) in
$\setR^{M_n}$. Then, theorem~\ref{thm:lower_bound} follows from
lemma~\ref{lem:H_n_deviation} if we prove that on $\mrm H_n$,
\begin{equation*}
  \inf_{\wh \tta} \int_{\Theta} \Prob_{\tta}^n \big\{
  \max_{j \in \mc J_n} |\wh \tta_j - \tta_j | \geq 1 - \von | \mf X_n
  \big\} \mc B(d\tta) \geq 1 - o(1),
\end{equation*}
or equivalently, that on $\mrm H_n$
\begin{equation}
  \label{eq:lower_bound_sup_argument}
  \sup_{\wh \tta} \int_{\Theta} \Prob_{\tta}^n \big\{
  \max_{j \in \mc J_n} |\wh \tta_j - \tta_j | < 1 - \von | \mf
  X_n \big\} \mc B(d\tta) = o(1).
\end{equation}
To prove~\eqref{eq:lower_bound_sup_argument}, we choose
\begin{equation*}
  \Theta = \Theta_{\von}^{M_n}, \quad
  \Theta_{\von} =  \{ -(1 - \von), 1 - \von \}, \quad \mc B =
  \bigotimes_{j \in \mc J_n} b_{\von}, \quad b_{\von} = \frac{1}{2}
  \big( \delta_{-(1 - \von)} + \delta_{1 - \von} \big),
\end{equation*}
where $\delta$ stands for the Dirac mass. Note that using
lemma~\ref{lem:lemma_y_j_sufficiency}, the left hand side
of~\eqref{eq:lower_bound_sup_argument} is smaller than
\begin{align*}
  \int \frac{\prod_{i=1}^n g_{\sigma}(Y_i)}{\prod_{j \in \mc J_n}
    g_{v_j}(y_j)} \Big( \prod_{j \in \mc J_n} \sup_{\wh \tta_j \in
    \setR} \int_{\Theta_{\von}} \ind{ | \wh \tta_j - \tta_j| < 1 -
    \von} g_{v_j} (y_j - \tta_j) db_{\von}(\tta_j) \Big) dY_1 \ldots d
  Y_n,
\end{align*}
and an easy argument shows that
\begin{equation*}
  \wh \tta_j = (1 - \von) \ind{y_j \geq 0} - (1 - \von) \ind{y_j < 0}
\end{equation*}
are strategies attaining the maximum. Thus, it suffices to prove the
lower bound among estimators $\wh \tta$ with coordinates $\wh \tta_j
\in \Theta_{\von}$ and measurable with respect to $y_j$ only. Since
the $y_j$ are independent with distribution density $g_{v_j}(\cdot -
\tta_j)$, the left hand side of~\eqref{eq:lower_bound_sup_argument} is
smaller than
\begin{multline*}
  \prod_{j \in \mc J_n} \max_{\wh \tta_j \in \Theta_{\von}}
  \int_{\Theta_{\von}} \int_{\setR} \ind{ | \wh \tta_j(u_j) - \tta_j|
    < 1 - \von} \, g_{v_j}(u_j - \tta_j) du_j \, db_{\von}
  (\tta_j) \\
  = \prod_{j \in \mc J_n} \Big( 1 - \inf_{\wh \tta_j \in
    \Theta_{\von}} \int_{\Theta_{\von}} \int_{\setR} \ind{| \wh
    \tta_j(u) - \tta_j | \geq 1 - \von} \, g_{v_j}(u - \tta_j) du \, d
  b_{\von}(\tta_j) \Big),
\end{multline*}
and if $\Phi(x) = \int_{-\infty}^x g_1(t)dt$ and $D_1$ is a positive
constant,
\begin{align*}
  \inf_{\wh \tta_j \in \Theta_{\von}} &\int_{\Theta_{\von}}
  \int_{\setR} \ind{| \wh \tta_j(u) - \tta_j | \geq 1 - \von} \,
  g_{v_j}(u - \tta_j) du \, db_{\von} \\
  &\geq \inf_{\wh \tta_j \in \Theta_{\von}} \frac{1}{2} \int_{\setR}
  \big( \ind{\wh \tta_j \geq 0} + \ind{\wh \tta_j < 0} \big)
  g_{v_j}(u - (1 - \von)) \wedge g_{v_j}(u + (1 - \von)) du \\
  &= \frac{1}{v_j} \int_{-\infty}^{0} g_{1}\Big(\frac{y - (1 -
    \von)}{v_j}\Big) du \\
  &= \Phi\Big( -\frac{1 - \von}{v_j} \Big) \geq \frac{D_1}{\sqrt{\log
      n}} n^{-(1 - \von)^2 (1 + \von) / (2s + 1)},
\end{align*}
where we used lemma~\ref{lem:lemma_y_j_sufficiency} and the fact that
for $x > 0$, $\Phi(-x) = \frac{(1 + o(1))\exp(-x^2 / 2) }{x
  \sqrt{2\pi}}$. It follows that the left hand side
of~\eqref{eq:lower_bound_sup_argument} is smaller than
\begin{multline*}
  \Big( 1 - \frac{D_1}{\sqrt{\log n}} n^{-(1 - \von)^2 (1 + \von) /
    (2s + 1)} \Big)^{M_n} \\
  \leq \exp\Big( |I_n| \, \Xi_n^{-1} \log\big( 1 - D_1 n^{-(1 -
    \von)^2 (1 + \von) / (2s + 1)} (\log n)^{-1 / 2} \big) \Big),
\end{multline*}
and if $D_2$ is a positive constant,
\begin{multline*}
  |I_n| \, \Xi_n^{-1} n^{-(1 - \von)^2 (1 + \von) / (2s + 1)} (\log
  n)^{-1/2} \\ = D_2 |I_n| \, n^{\von / (2s + 1)} \times n^{\von^2 (1
    - \von) / (2s + 1)} (\log n)^{-1/2 - 1 / (2s + 1)} \raro +\infty
\end{multline*}
as $n \raro +\infty$, since $|I_n| n^{\von / (2s + 1)} \raro +\infty$,
thus the theorem. $\hfill \qed$

\appendix

\section{Well known facts on optimal recovery}
\label{sec:optimal_recovery}

\subsection{Explicit values}
\label{sec:explicit_values}

\begin{figure}[htbp]
  \begin{center}
    \includegraphics[width = 12cm]{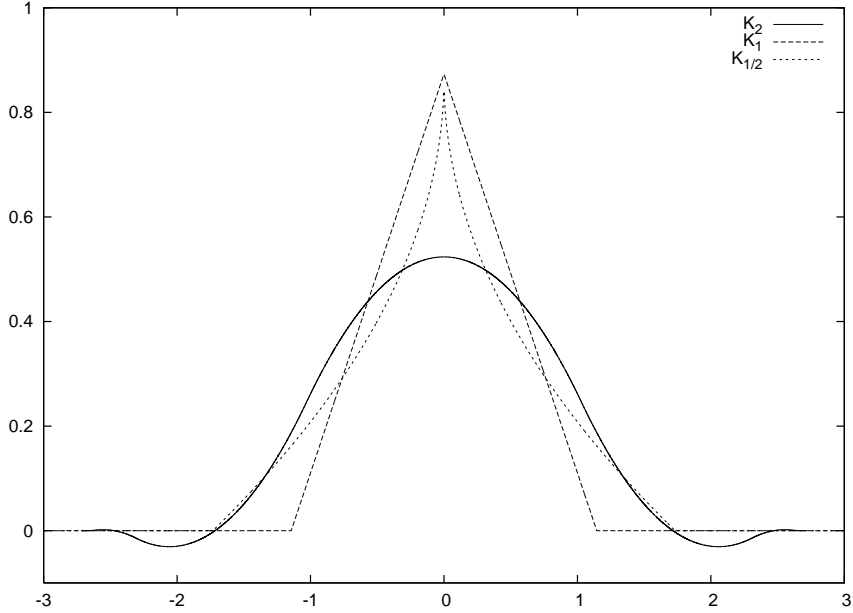}
  \end{center}
  \caption{Optimal recovery kernels $K_s$ for $s = 1/2$, $s = 1$ and
    $s = 2$.}
  \label{fig_kernels}
\end{figure}

To our knowledge, the function $\varphi_s$ is only known for $s \in
(0,1] \cup \{ 2 \}$. We recall that the optimal recovery kernel is
defined by
\begin{equation*}
  K_s = \frac{\varphi_s}{\int_{\setR} \varphi_s},
\end{equation*}
where $\varphi_s$ is given by~\eqref{eq:OR_optimal_function}. The
kernel $K_s$ for $s \in (0, 1]$ was found by \cite{korostelev93} and
\cite{fuller60} for $s = 2$. See also \cite{leonov97, leonov99},
\cite{lepski_tsybakov00} and \cite{bertin_phd}. When $s \in (0, 1]$,
\begin{equation*}
  K_s(t) = \frac{s+1}{2s} \varphi_s^{-1 / s}(0) \big(1 -
  \varphi_s^{-1}(0)|t|^s \big)_+,
\end{equation*}
where $x_+ = \max(0, x)$, and
\begin{equation*}
  \varphi_s(0) = \Big( \frac{(2s + 1)(s + 1)}{4 s^2} \Big)^{s / (2s +
    1)}.
\end{equation*}
When $s = 2$, we have
\begin{equation*}
  \varphi_s(t) = \tta^{-2/5} g_2(\tta^{-2/5} t),
\end{equation*}
where for $t \geq 0$
\begin{align*}
  g_2(t) &= \sum_{j \geq 0} \big( (-1)^j q^j + \frac{1}{2} (-1)^{j+1}
  (t - t_{2j})^2 \big) \ind{t \in [t_{2j - 1}, t_{2j + 1}]}, \\
  q &= \frac{1}{16} \Big(3 + \sqrt{33} - \sqrt{ 26 + 6 \sqrt{33}} \,\,
  \Big)^2, \\
  \tta &= \frac{2(23 q^2 - 14 q + 23) \sqrt{1 + q}}{30 (1 - q^{5/2})},
\end{align*}
and $t_{-1} = t_{0} = 0$, $t_1 = \sqrt{1 + q}$ and for any $j \in
\setN - \{ 0 \}$, $t_{2j} = 2 \sqrt{1 + q} \sum_{i = 0}^{j-1}
q^{i/2}$, $t_{2j + 1} = t_{2j} + q^{j/2} \sqrt{1 + q}$. Note that
$\varphi_2$ is piecewise quadratic and infinitely oscillating around
$0$ at the boundaries of its support. For these values of $s$,
\begin{equation*}
  P = P_s =
  \begin{cases}
    \, \displaystyle \Big( \frac{s + 1}{2 s^2} \Big)^{s / (2s + 1)}
    &\text{ when } s \in (0, 1], \\[0.3cm]
    \, \displaystyle \Big( \frac{2}{5} \Big)^{2/5} \tta^{-2/5} &\text{
      when } s = 2.
  \end{cases}
\end{equation*}
In figure~\ref{fig_kernels} we give an illustration of the kernel
$K_s$ for $s = 1/2$, $s = 1$ and $s = 2$.

\subsection{Optimal recovery}

The next results are well known and can be found in \cite{donoho94},
\cite{leonov97, leonov99}, \cite{lepski_tsybakov00} and
\cite{bertin_phd}. The problem consists in recovering $f$ from
\begin{equation}
  \label{eq:optimal_recovery_model}
  y(t) = f(t) + \von z(t), \quad t \in \setR,
\end{equation}
where $\von > 0$, $z$ is an unknown deterministic function such that
$\norm{z}_2 \leq 1$ and $f \in C(s, L; \setR) \eqdef \Sigma(s, L;
\setR) \cap \mbb L^2(\setR)$. This problem is well known, and the link
between this problem and the statistical estimation in sup norm in the
white noise model
\begin{equation*}
  d Y_t^{\von} = f(t) dt + \von dW_t, \quad t \in \setR,
\end{equation*}
was made by \cite{donoho94}, see also \cite{leonov99}. The minimax
error for the problem of optimal recovery of $f$ at $0$ in the
model~\eqref{eq:optimal_recovery_model} is defined by
\begin{equation*}
  E_s(\von, L) \eqdef \inf_{T} \sup_{\substack{ f \in C(s, L; \setR)
      \\ \norm{f - y}_2 \leq \von}} | T(y) - f(0) |,
\end{equation*}
where $\inf_T$ is taken among all continuous and linear forms on $\mbb
L^2(\setR)$. We know from \cite{miccheli_rivlin77}, \cite{arestov90}
that
\begin{align*}
  E_s(\von, L) &= \inf_{K \in \mbb L^2(\setR)} \Big( \sup_{f \in C(s,
      L; \setR) } \Big| \int K(t) (f(t) - f(0)) \Big| + \von
      \norm{K}_2 \Big) \\ &= \sup_{\substack{f \in \Sigma(s, L; \setR)
      \\ \norm{f}_2 \leq \von}} f(0).
\end{align*}
Note that $\varphi_s$ satisfies $\varphi_s(0) = E_s(1, 1)$. For any $s
> 0$, we know from \cite{leonov97} that $\varphi_s$ is well defined
and unique, that it is even and compactly supported and that
$\norm{\varphi_s}_2 = 1$. A renormalisation argument from
\cite{donoho94} shows that
\begin{equation*}
  E_s(\von, L) = E_s(1, 1) L^{1/(2s+1)} \von^{2s / (2s+1)},
\end{equation*}
thus it suffices to know $E_s(1, 1)$. If we define
\begin{equation}
  \label{eq:B_s_L_def}
  \mc B(s, L) \eqdef \sup_{f \in C(s, L; \setR)} \Big| \int K_s(t) (
  f(t) - f(0) ) \Big|,
\end{equation}
we have the decomposition
\begin{equation*}
  E_s(1, 1) = \mc B(s, 1) + \norm{K}_2,
\end{equation*}
and in particular if $P$ is given by~\eqref{eq:P_s_definition} and
$c_s$ by~\eqref{eq:c_s_def} we have
\begin{equation}
  \label{eq:P_s_decomposition}
  P = L c_s^s \big( \mc B(s, 1) + \norm{K}_2 \big).
\end{equation}


\footnotesize

\bibliographystyle{ims}

\bibliography{biblio}



\end{document}

%% file: notation.tex
\newtheorem{proposition}{Proposition}%
\newtheorem{lemma}{Lemma}%
\newtheorem{theorem}{Theorem}

\theoremstyle{definition}

\newcounter{assumption}
\newtheorem{assumption}[assumption]{Assumption}

\renewcommand{\leq}{\leqslant}
\renewcommand{\geq}{\geqslant}

\DeclareMathOperator*{\argmax}{argmax}
\DeclareMathOperator*{\argmin}{argmin}
\DeclareMathOperator*{\Span}{Span}

\newcounter{hyp}

\newcommand{\mc}{\mathcal}
\newcommand{\mb}{\mathbf}
\newcommand{\mbb}{\mathbb}
\newcommand{\mrm}{\mathrm}
\newcommand{\mf}{\mathfrak}

\newcommand{\tup}[1]{\textup{#1}}

\newcommand{\wh}{\widehat}

\newcommand{\wt}{\widetilde}

\newcommand{\kpa}{\kappa}

\newcommand{\von}{\varepsilon}

\newcommand{\tta}{\theta}

\newcommand{\Lba}{\Lambda}

\newcommand{\lba}{\lambda}

\newcommand{\raro}{\rightarrow}

\newcommand{\sumin}{\sum_{i=1}^{n}}

\newcommand{\setR}{\mathbb{R}}

\newcommand{\setN}{\mathbb{N}}

\newcommand{\eqdef}{\triangleq}

\newcommand{\ind}[1]{\mathbf{1}_{#1}}

\newcommand{\prodsca}[2]{\langle #1 \, ,\, #2 \rangle}

\newcommand{\norm}[1]{\| #1 \|}

\newcommand{\norminfty}[1]{\| #1 \|_{\infty}}

\newcommand{\supp}[1]{\text{Supp }#1}

\newcommand{\ppint}[1]{\lfloor #1 \rfloor} 

\newcommand{\Prob}{\mbb{P}}

\newcommand{\Pm}{\Prob_{\mu}^{n}}
\newcommand{\Pfm}{\Prob_{f,\mu}^{n}}
\newcommand{\Efm}{\mbb E_{\,f, \mu}^{\,n}}

\newcommand{\Em}{\mbb E_{\mu}^{n}}

\newcommand{\Var}{\mbb V \text{ar}}





%% file: article3.bbl
\begin{thebibliography}{33}
\expandafter\ifx\csname natexlab\endcsname\relax\def\natexlab#1{#1}\fi
\expandafter\ifx\csname url\endcsname\relax
  \def\url#1{\texttt{#1}}\fi
\expandafter\ifx\csname urlprefix\endcsname\relax\def\urlprefix{URL }\fi
\providecommand{\eprint}[2][]{\url{#2}}

\bibitem[{Antoniadis et~al.(1997)Antoniadis, Gregoire and
  Vial}]{antoniadis_et_al97}
\textsc{Antoniadis, A.}, \textsc{Gregoire, G.} and \textsc{Vial, P.} (1997).
\newblock Random design wavelet curve smoothing.
\newblock \textit{Statistics and Probability Letters}, \textbf{35} 225--232.

\bibitem[{Arestov(1990)}]{arestov90}
\textsc{Arestov, V.~V.} (1990).
\newblock Optimal recovery of operators and related problems.
\newblock \textit{Proc. Steklov Inst. Math.}, \textbf{4} 1--20.

\bibitem[{Bertin(2004{\natexlab{a}})}]{bertin04}
\textsc{Bertin, K.} (2004{\natexlab{a}}).
\newblock Asymptotically exact minimax estimation in sup-norm for anisotropic
  h\"older classes.
\newblock \textit{Bernoulli}, \textbf{10} 873--888.

\bibitem[{Bertin(2004{\natexlab{b}})}]{bertin_phd}
\textsc{Bertin, K.} (2004{\natexlab{b}}).
\newblock \textit{Estimation asymptotiquement exacte en norme sup de fonctions
  multidimensionnelles}.
\newblock Ph.D. thesis, Universit\'e Paris 6.

\bibitem[{Bertin(2004{\natexlab{c}})}]{bertin02}
\textsc{Bertin, K.} (2004{\natexlab{c}}).
\newblock Minimax exact constant in sup-norm for nonparametric regression with
  random design.
\newblock \textit{J. Statist. Plann. Inference}, \textbf{123} 225--242.

\bibitem[{Brown and Cai(1998)}]{cai_brown98}
\textsc{Brown, L.} and \textsc{Cai, T.} (1998).
\newblock Wavelet shrinkage for nonequispaced samples.
\newblock \textit{The Annals of Statistics}, \textbf{26} 1783--1799.

\bibitem[{Brown et~al.(2002)Brown, Cai, Low and Zhang}]{brown_cai_low_zhang02}
\textsc{Brown, L.~D.}, \textsc{Cai, T.}, \textsc{Low, M.~G.} and \textsc{Zhang,
  C.-H.} (2002).
\newblock Asymptotic equivalence theory for nonparametric regression with
  random design.
\newblock \textit{The Annals of Statistics}, \textbf{30} 688 -- 707.

\bibitem[{Brown and Low(1996)}]{brown_low96}
\textsc{Brown, L.~D.} and \textsc{Low, M.~G.} (1996).
\newblock Asymptotic equivalence of nonparametric regression and white noise.
\newblock \textit{The Annals of Statistics}, \textbf{24} 2384--2398.

\bibitem[{Cai and Low(2004{\natexlab{a}})}]{cai_low04a}
\textsc{Cai, T.~T.} and \textsc{Low, M.~G.} (2004{\natexlab{a}}).
\newblock An adaptation theory for nonparametric confidence intervals.
\newblock \textit{The Annals of Statistics}, \textbf{32} 1805--1840.

\bibitem[{Cai and Low(2004{\natexlab{b}})}]{cai_low04b}
\textsc{Cai, T.~T.} and \textsc{Low, M.~G.} (2004{\natexlab{b}}).
\newblock Adaptive confidence balls.
\newblock \textit{The Annals of Statistics}.
\newblock To appear.

\bibitem[{Donoho(1994)}]{donoho94}
\textsc{Donoho, D.~L.} (1994).
\newblock Asymptotic minimax risk for sup-norm loss: Solution via optimal
  recovery.
\newblock \textit{Probability Theory and Related Fields}, \textbf{99} 145--170.

\bibitem[{Donoho and Low(1992)}]{donoho_low_92}
\textsc{Donoho, D.~L.} and \textsc{Low, M.~G.} (1992).
\newblock Renormalization exponents and optimal pointwise rates of convergence.
\newblock \textit{The Annals of Statistics}, \textbf{20} 944--970.

\bibitem[{Fan and Gijbels(1995)}]{fan_gijbels95}
\textsc{Fan, J.} and \textsc{Gijbels, I.} (1995).
\newblock Data-driven bandwidth selection in local polynomial fitting: variable
  bandwidth and spatial adaptation.
\newblock \textit{Journal of the Royal Statistical Society. Series B.
  Methodological}, \textbf{57} 371--394.

\bibitem[{Fan and Gijbels(1996)}]{fan_gijbels96}
\textsc{Fan, J.} and \textsc{Gijbels, I.} (1996).
\newblock \textit{Local polynomial modelling and its applications}.
\newblock Monographs on Statistics and Applied Probability, Chapman \& Hall,
  London.

\bibitem[{Fuller(1961)}]{fuller60}
\textsc{Fuller, A.~T.} (1961).
\newblock Relay control systems optimized for various performance criteria,.
\newblock \textit{Automatic and remote control}, \textbf{1}.

\bibitem[{Ga\"iffas(2004)}]{gaiffas04a}
\textsc{Ga\"iffas, S.} (2004).
\newblock Convergence rates for pointwise curve estimation with a degenerate
  design.
\newblock \textit{Mathematical Methods of Statistics}.
\newblock To appear, available at http://hal.ccsd.cnrs.fr/ccsd-00003086/en/.

\bibitem[{Hoffmann and Lepski(2002)}]{hoffmann_lepski02}
\textsc{Hoffmann, M.} and \textsc{Lepski, O.~V.} (2002).
\newblock Random rates in anisotropic regression.
\newblock \textit{The Annals of Statistics}, \textbf{30} 325--396.

\bibitem[{Korostelev and Nussbaum(1999)}]{korostelev_nussbaum_99}
\textsc{Korostelev, A.} and \textsc{Nussbaum, M.} (1999).
\newblock The asymptotic minimax constant for sup-norm loss in nonparametric
  density estimation.
\newblock \textit{Bernoulli}, \textbf{5} 1099--1118.

\bibitem[{Korostelev(1993)}]{korostelev93}
\textsc{Korostelev, V.} (1993).
\newblock An asymptotically minimax regression estimator in the uniform norm up
  to exact contant.
\newblock \textit{Theory of Probability and its Applications}, \textbf{38}
  737--743.

\bibitem[{Korostelev and Tsybakov(1993)}]{korostelev_tsybakov93}
\textsc{Korostelev, V.} and \textsc{Tsybakov, A.} (1993).
\newblock \textit{Minimax theory of image reconstruction}.
\newblock Springer-Verlag, New York.

\bibitem[{Ledoux and Talagrand(1991)}]{ledoux_talagrand91}
\textsc{Ledoux, M.} and \textsc{Talagrand, M.} (1991).
\newblock \textit{Probability in {B}anach spaces}, vol.~23 of
  \textit{Ergebnisse der Mathematik und ihrer Grenzgebiete (3) [Results in
  Mathematics and Related Areas (3)]}.
\newblock Springer-Verlag, Berlin.
\newblock Isoperimetry and processes.

\bibitem[{Leonov(1997)}]{leonov97}
\textsc{Leonov, S.} (1997).
\newblock On the solution of an optimal recovery problem and its applications
  in nonparametric regression.
\newblock \textit{Mathematical Methods of Statistics}, \textbf{6} 476--490.

\bibitem[{Leonov(1999)}]{leonov99}
\textsc{Leonov, S.} (1999).
\newblock Remarks on extremal problems in nonparametric curve estimation.
\newblock \textit{Statistics and Probability Letters}, \textbf{43} 169--178.

\bibitem[{Lepski and Tsybakov(2000)}]{lepski_tsybakov00}
\textsc{Lepski, O.~V.} and \textsc{Tsybakov, A.~B.} (2000).
\newblock Asymptotically exact nonparametric hypothesis testing in sup-norm and
  at a fixed point.
\newblock \textit{Probability Theory and Related Fields}, \textbf{117} 17--48.

\bibitem[{Low(1997)}]{low97}
\textsc{Low, M.~G.} (1997).
\newblock On nonparametric confidence intervals.
\newblock \textit{The Annals of Statistics}, \textbf{25} 2547--2554.

\bibitem[{Maxim(2003)}]{voichitaphd}
\textsc{Maxim, V.} (2003).
\newblock \textit{Restauration de signaux bruit\'es sur des plans d'experience
  al\'eatoires}.
\newblock Ph.D. thesis, Universit\'e Joseph Fourier, Grenoble 1.

\bibitem[{Micchelli and Rivlin(1977)}]{miccheli_rivlin77}
\textsc{Micchelli, C.~A.} and \textsc{Rivlin, T.~J.} (1977).
\newblock A survey of optimal recovery.
\newblock \textit{Optimal estimation in approximation theory} 1 -- 54.

\bibitem[{Picard and Tribouley(2000)}]{picard_tribouley00}
\textsc{Picard, D.} and \textsc{Tribouley, K.} (2000).
\newblock Adaptive confidence interval for pointwise curve estimation.
\newblock \textit{The Annals of Statistics}, \textbf{28} 298--335.

\bibitem[{Pinsker(1980)}]{pinsker80}
\textsc{Pinsker, M.~S.} (1980).
\newblock Optimal filtration of functions from ${L}_2$ in {G}aussian noise.
\newblock \textit{Problems of Information Transmission}, \textbf{16} 52--68.

\bibitem[{Spokoiny(1998)}]{spok98}
\textsc{Spokoiny, V.~G.} (1998).
\newblock Estimation of a function with discontinuities via local polynomial
  fit with an adaptive window choice.
\newblock \textit{The Annals of Statistics}, \textbf{26} 1356--1378.

\bibitem[{Stone(1982)}]{stone82}
\textsc{Stone, C.~J.} (1982).
\newblock Optimal global rates of convergence for nonparametric regression.
\newblock \textit{The Annals of Statistics}, \textbf{10} 1040--1053.

\bibitem[{Tsybakov(2003)}]{tsybakov03}
\textsc{Tsybakov, A.} (2003).
\newblock \textit{Introduction à l'estimation non-paramétrique}.
\newblock Springer.

\bibitem[{Wong and Zheng(2002)}]{wong_zheng02}
\textsc{Wong, M.-Y.} and \textsc{Zheng, Z.} (2002).
\newblock Wavelet threshold estimation of a regression function with random
  design.
\newblock \textbf{80} 256--284.

\end{thebibliography}
